\documentclass{article}

\usepackage{PRIMEarxiv}

\usepackage[utf8]{inputenc} 
\usepackage[T1]{fontenc}    
\usepackage{hyperref}       
\usepackage{url}            
\usepackage{booktabs}       
\usepackage{amsfonts}       
\usepackage{nicefrac}       
\usepackage{microtype}      
\usepackage{lipsum}
\usepackage{fancyhdr}       
\usepackage{graphicx}       
\graphicspath{{media/}}     
\usepackage{amsmath}
\usepackage{multirow}
\usepackage{algorithm}
\usepackage{algpseudocode}  

\title{A quantum algorithm for the Kalman filter using block encoding

}

\author{
  Hao Shi\\
  College of Intelligence Science \\
  National University of Defense Technology \\
  Changsha, China\\
  \texttt{shihao\_hill@163.com} \\
   \And
  Guofeng Zhang \\
  Department of Applied Mathematics \\
  Hong Kong Polytechnic University \\
  Hong Kong, China\\
  \texttt{guofeng.zhang@polyu.edu.hk} \\
  \And
  Ming Zhang* \\
  College of Intelligence Science \\
  National University of Defense Technology \\
  Changsha, China\\
  \texttt{zhangming@nudt.edu.cn} \\
}

\begin{document}
\maketitle

\begin{abstract}
Quantum algorithms offer significant speed-ups over their classical counterparts in various applications. In this paper, we develop quantum algorithms for the Kalman filter widely used in classical control engineering using the block encoding method. The entire calculation process is achieved by performing matrix operations on Hamiltonians based on the block encoding framework, including addition, multiplication, and inversion, which can be completed in a unified framework compared to previous quantum algorithms for solving control problems. We demonstrate that the quantum algorithm exponentially accelerates the computation of the Kalman filter compared to traditional methods. The time complexity can be reduced from $O(n^3)$ to $O(\kappa poly\log(n/\epsilon)\log(1/\epsilon'))$, where $n$ represents the matrix dimension, $\kappa$ denotes the condition number for the matrix to be inverted, $\epsilon$ indicates desired precision in block encoding, $\epsilon'$ signifies desired precision in matrix inversion. This paper provides a comprehensive quantum solution for implementing the Kalman filter and serves as an attempt to broaden the scope of quantum computation applications. Finally, we present an illustrative example implemented in Qiskit (a Python-based open-source toolkit) as a proof-of-concept.
\end{abstract}


\section{Introduction}
Quantum computing is an emerging and rapidly growing paradigm that has the potential to offer significant computational advantages over classical computing, by leveraging quantum-mechanical principles such as quantum superposition and entanglement\cite{nielsen_chuang_2000,Gill2022}. This computational advantage of quantum computing will facilitate the resolution of complex computational problems in various application domains. As is well known, classical computers store information in bits (`0' or `1'), whereas their equivalent in quantum computers are qubits, represented as $|0\rangle$ or $|1\rangle$ or any superposition of $|0\rangle$ and $|1\rangle$. By utilizing Hilbert space as a quantum computing space, where $n$ qubits can result in a superposition state of $2^n$ possible outcomes, quantum computers possess an inherent advantage over classical computers when it comes to solving large-scale problems. Over the past few decades, numerous quantum algorithms have been proposed, two well-known algorithms are Shor's algorithm for factoring\cite{Shor1994} and Grover's algorithm for database searching\cite{Grover1996}, which are often employed to illustrate the advantages of quantum computing. The core steps of both algorithms use the subroutines of quantum phase estimation \cite{nielsen_chuang_2000} and amplitude amplification \cite{Brassard2000}. In 2009, the Harrow-Hassidim-Lloyd (HHL) algorithm was proposed for linear systems\cite{Harrow2009}, which has the potential of solving a large number of mathematical problems, such as linear differential equations\cite{Berry2012}, least-square curve fitting\cite{Wiebe2012} and matrix inversions\cite{TaShma2013}. Many subsequent quantum algorithms implemented using these methods as subroutines achieve exponential speedup $O(log(n))$ or quadratic speedup $O(\sqrt{n})$, relative to their classical counterparts \cite{Biamonte2017}.

As linear system problems arise ubiquitously across all scientific and engineering disciplines, one may expect to explore more quantum algorithms. Notably, Wossnig et al.\cite{Wossnig2017} proposed a quantum algorithm for dense matrices based on HHL without considering matrix sparsity. Gily{\'e}n et al. introduced a technique called ``block encoding", which involves embedding a general matrix into a larger unitary operator implemented by the quantum system using blocks of unitaries. They also established rules for addition, subtraction and multiplication of block-encoded matrices, thereby leading to exponential speed-ups with respect to matrix dimension\cite{Gilyen2019}. Subsequently, a hybrid quantum-classical algorithm called variational quantum linear solver (VQLS) is presented by Bravo-Prieto et al.\cite{Bravo-Prieto2019}, which used fewer qubits compared with HHL. Shao\cite{Shao2018} explored three approaches for solving matrix multiplication operations and discussed the calculation problem from the perspective of input/output data types (quantum/classical). Li et al.\cite{Li2021} proposed a quantum matrix multiplier that can serve as an intermediate module for more complex quantum algorithms. Qi et al.\cite{Qi2022} developed quantum algorithms for matrix operations using the Sender-Receiver model. Zhao et al.\cite{Zhao2021} provided Trotter-based subroutines for elementary matrix operations through the matrix embedding formula. These aforementioned quantum algorithms hold promise in terms of offering speed-ups over their classical counterparts when applied to specific problems.

Matrix operations play a crucial role in machine learning and data processing, and quantum computing holds immense potential for enhancing matrix arithmetic, as mentioned above. In the field of control theory and engineering, researchers have explored the application of quantum computing to solve linear system problems encountered in classical control systems. Sun et al.\cite{Sun2017} proposed a quantum algorithm based on HHL for solving Lyapunov equations, demonstrating exponential acceleration of optimal control protocols using quantum algorithms\cite{Sun2016}. Li et al.\cite{Li2020} proposed a quantum algorithm that enhances the efficiency of solving state estimation problems. Furthermore, Li et al.\cite{Li2022} devised a quantum linear system algorithm for general matrices in system identification. All these aforementioned quantum algorithms exhibit accelerated performance compared to their classical counterparts. In this paper, we intent to develop a  block-encoding-based quantum version of the classical Kalman filter, showcasing its exponential improvement.

Being the best Bayesian estimator of optimal state estimation for linear systems with Gaussian uncertainty\cite{Kim2018}, the Kalman filter is widely employed in various fields such as in UAV, robotics, communication, navigation\cite{Welch1995}. It is a recursive estimator that combines mathematical models and output measurements from a physical system to accurately estimate its state. Solving state estimation involves solving a series of linear equations, which can be efficiently addressed by quantum algorithms. Therefore, the proposed quantum algorithm for the Kalman filter exponentially reduce time complexity from $O(n^3)$ to $O(\kappa poly\log(n/\epsilon)\log(\kappa/\epsilon'))$.

This paper is organized as follows: Section 2 presents the proposed block encoding schemes and matrix operations based on block encoding. Subsequently, in Section 3, our quantum algorithm and the corresponding quantum circuits for the Kalman filter are introduced. In Section 4, a numerical example is provided to demonstrate the feasibility of our approach using Qiskit\cite{Qiskit2023}. Finally, Section 5 concludes the paper.

\section{Block encoding}
\begin{figure}[!t]
	\centering
	\includegraphics[scale=0.3]{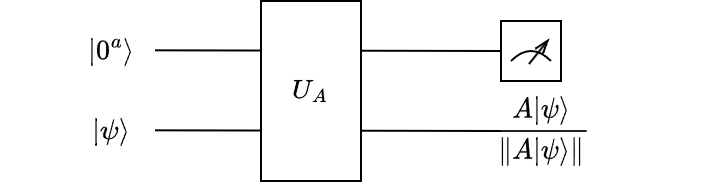}
	\caption{Circuit for block encoding of $A$ using $a$ ancilla qubits.}
	\label{fig:block_encoding_01}
\end{figure}

Block encoding is a general way of transforming a matrix to a unitary operator. An illustrative example of block encoding is as follows: if an $s$-qubit matrix $A$ is embedded in an $(s+a)$-qubit unitary matrix $U_A$,
\begin{equation}
    U_A=
	\begin{bmatrix}
		A/\alpha & * & \cdots & * \\
		* & * & \cdots & * \\
		\vdots & & \vdots & \\
		* & * & \cdots & * 
	\end{bmatrix},
\end{equation} 
it represents a submatrix with the normalization factor $\alpha$ as the upper-left block of $U_A$, where $*$ denotes irrelevant $s$-qubit matrices and $U_A$ consists of $2^a$ block rows(columns). Now suppose that $|\psi\rangle$ is any quantum state consisting of $s$ qubits, and $|0^a\rangle$ represents $a$ ancilla qubits. By measuring the ancilla qubits, we can obtain outcome $A$ if it returns 0. The circuit is presented in Figure \ref{fig:block_encoding_01}. The equation for solving the matrix $A$ is given by:
\begin{align} 
  A=\alpha(\langle{0^a}|\otimes I_s)U_A(|0^a\rangle \otimes I_s)
\end{align}

\textbf{Theorem 1:} Block encoding\cite{Gilyen2019,Lin2022}. {Given an $s$-qubit matrix $A\in\mathbb{C}^{n\times n}$, if there exist real positive values $\alpha$,$\epsilon$ $\in \mathbb{R}_+$, a natural number $a \in \mathbb{N}$, and a unitary matrix $U_A$ of $(s+a)$ qubits, such that
	\begin{align}
		\Vert A-\alpha(\langle{0^a}|\otimes I_s)U_A(|0^a\rangle \otimes I_s) \Vert \le \epsilon,
	\end{align}
then $U_A$ is called an $(\alpha,a,\epsilon)$-block-encoding of $A$. In the case where the block encoding is exact with $\epsilon=0$, $U_A$ is termed an $(\alpha,a)$-block-encoding of A. The set of all $(\alpha,a,\epsilon)$-block-encoding of $A$ is denoted by $BE_{\alpha,a}(A,\epsilon)$, and we define $BE_{\alpha,a}(A)=BE(A,0)$. 

If $A\in\mathbb{C}^{m\times n}$, where $m,n\leq 2^s$, we can define $A_e\in \mathbb{C}^{2^s\times 2^s}$ such that $A$ occupies the upper-left of matrix $A_e$ while all other elements are set to zero \cite{Gilyen2019}.

\subsection{The construction of block encoding}

There exist many specific methods to find $U_A$ to block encode $A$\cite{Lin2022,Camps2024,Camps2022FABLE,Clader2022,Wan2019}. In this paper, building upon the theory in the references\cite{Gilyen2019,Lin2022}, we derive a method for constructing $U_A$ in quantum-accessible data structure introduced in \cite{Kerenidis2016} to efficiently generate quantum states.

\textbf{Theorem 2:} Block encoding of matrix in quantum data structure. Let $A\in \mathbb{C}^{m\times n}$ be an $s$-qubit matrix with $A_{ij}\in \mathbb{C}$ being the entry on the $i$-row and the $j$-column, let $m,n\leq 2^s$. When $A$ is stored in a quantum data structure, then we always can implement unitaries $U_R$ and $U_L$ in time $O(polylog(mn/\epsilon))$ in Theorem 1 according to \cite{Chakraborty2018}. 

\begin{align}
	U_L: |0^s\rangle|j\rangle\rightarrow \frac{1}{\Vert A\Vert_F}\sum_{i=1}^{m}\Vert A_{i,\cdot}\Vert |i\rangle |j\rangle \\
	U_R: |i\rangle|0^{s}\rangle \rightarrow |i\rangle\frac{1}{\Vert A_{i,\cdot}\Vert}\sum_{j=1}^{n}A_{i,j}|j\rangle
\end{align}
where $\Vert A_{i,\cdot} \Vert$ is the $l_2$ norm of the $i$-th row of $A$. So we can verify that $U_L^\dagger U_R$ is a $(\Vert A\Vert_F,s,\epsilon)$-block-encoding of $A$, where $\Vert A\Vert_F=\sqrt{\sum_{i,j}|A_{i,j}|^2}$ is the Frobenius norm of $A$. The construction circuit can be seen in Figure \ref{fig:Construction_of_BE}. $|0^s\rangle$ denotes the ancilla qubits used to construct $U_A$.

\begin{figure}[!h]
	\centering
	\includegraphics[scale=0.3]{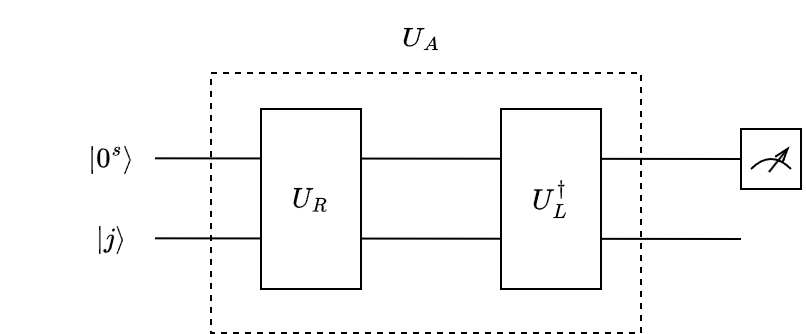}
	\renewcommand{\figurename}{Figure} 
	\caption{Circuit construction for block encoding of $A$ using $s$ ancilla qubits.}
	\label{fig:Construction_of_BE}
\end{figure}

\textit{Proof.}

\begin{align}
	\langle i|\langle0^s| U_L^\dagger U_R |0^s\rangle|j\rangle &= (U_L|i \rangle|0^s\rangle)^\dagger U_R |0^s\rangle|j\rangle = (\frac{1}{\Vert A\Vert_F}\sum_{i=1}^{m}\Vert A_{i,\cdot}\Vert \langle i| \langle j|) (|i\rangle\frac{1}{\Vert A_{i,\cdot}\Vert}\sum_{j=1}^{n}A_{i,j}|j\rangle) \nonumber\\
	&=\frac{A_{i,j}}{\Vert A\Vert_F}
\end{align}

So we consider that a $BE_{\Vert A\Vert_F,s}(A,\epsilon)$ can be implemented in time $O(polylog(mn/\epsilon))$ described above.

\subsection{Addition and multiplication of block-encoded matrices}
If the block encoding of a matrix has been efficiently completed, we can perform matrix addition and multiplication operations using the following rules. The theorems referenced in \cite{Gilyen2019,Lin2022} have been modified and expanded to adjust the specific problems addressed in our study.

\textbf{Theorem 3:} Addition of two block-encoded matrices. {If $U_A$ represents an ($\alpha,a,\delta$)-block-encoding of an $s$-qubit operator $A$ and $U_B$ represents a ($\beta,b,\epsilon$)-block-encoding of an $s$-qubit operator $B$, $\alpha,\beta >0$. Let $T=\alpha U_A+\beta U_B$ be the linear combination of unitaries (LCU), note that the dimensions of $U_A,U_B$ may not be equal and the sum is simply the addition of the corresponding elements. Define $W=(V^\dagger\otimes I_s)U(V\otimes I_s)$, we can get that $W\in BE_{\Vert\gamma\Vert_1,max(a,b)+1}(A+B,\delta+\epsilon)$, where $\Vert\gamma\Vert_1=\alpha+\beta$.

\textit{proof.} Firstly	
\begin{align}
	U_A=\begin{bmatrix}
		A'/\alpha & * & \cdots & *\\
		* & * & \cdots & * \\
		\vdots & & \vdots &  \\
		* & * & \cdots & *
	\end{bmatrix},
	U_B=\begin{bmatrix}
	B'/\alpha & * & \cdots & *\\
	* & * & \cdots & * \\
	\vdots & & \vdots &  \\
	* & * & \cdots & *
\end{bmatrix}
	\label{equ:U}
\end{align}
and 
\begin{align}
	\Vert A-A'\Vert \leq \delta, \Vert B-B'\Vert \leq \epsilon.
\end{align}
\begin{align}
	T=\alpha U_A+\beta U_B=
	\begin{bmatrix}
		A'+B' & * & \cdots & *\\
		* & * & \cdots & * \\
		\vdots & & \vdots &  \\
		* & * & \cdots & *
	\end{bmatrix}.
\end{align}
even if the dimensions of $U_A$ and $U_B$ are not equal, where $T\in\mathbb{C}^{2^{\max(a,b)}\times 2^{\max(a,b)}}$. So
\begin{equation}
	\begin{aligned}
	\Vert A+B-(\langle 0^{max(a,b)}|\otimes I_s)T(|0^{max(a,b)}\rangle\otimes I_s)\Vert &=\Vert A+B-(A'+B')\Vert \\
	& \leq \Vert A-A'\Vert +\Vert B-B'\Vert \\
	& \leq \delta+\epsilon.
	\end{aligned}
\end{equation}
Secondly, let
\begin{align}
		U=|0\rangle\langle 0|\otimes U_A + |1\rangle\langle 1|\otimes U_B,
		\label{equ:U}
\end{align}
and
\begin{align}
	V=\frac{1}{\sqrt{\Vert\gamma\Vert_1}}
	\begin{bmatrix}
		\sqrt{\alpha} & * \\
		\sqrt{\beta} & * 
	\end{bmatrix},
	V^\dagger=\frac{1}{\sqrt{\Vert\gamma\Vert_1}}
	\begin{bmatrix}
		\sqrt{\alpha} & \sqrt{\beta} \\
		* & * 
	\end{bmatrix}.
	\label{equ:V}
\end{align}
According to Lemma 7.9 in \cite{Lin2022}, we can verify that $W\in BE_{\Vert\gamma\Vert_1,1}(T)$, where $\Vert\gamma\Vert_1=\alpha+\beta$. 

Thirdly, 
\begin{align}
	W=
	\begin{bmatrix}
		\frac{T}{\Vert\gamma\Vert_1} & * \\
		* & * 
	\end{bmatrix}
	\label{equ:W}
\end{align}
where $W\in \mathbb{C}^{2^{\max(a,b)+1}\times 2^{\max(a,b)+1}}$.

Lastly,
\begin{equation}
	\begin{aligned}
		& \Vert A+B- \Vert \gamma\Vert_1 (\langle0|\langle 0^{max(a,b)}| \otimes I_s)W(|0\rangle|0^{max(a,b)}\rangle \otimes I_s)\Vert \\
		=&  \Vert A+B- (\langle 0^{max(a,b)}| \otimes I_s)T(|0^{max(a,b)}\rangle \otimes I_s) \Vert \\
		 \leq & \delta+\epsilon.
		\label{equ:delta}
	\end{aligned}
\end{equation}

Therefore, $W$ serves as a $({\Vert\gamma\Vert_1,max\{a,b\}+1,\delta+\epsilon})$-block-encoding of $A+B$. The quantum circuit is illustrated in Figure \ref{fig:Addition01}.
\begin{figure}[!t]
	\centering
	\includegraphics[scale=0.3]{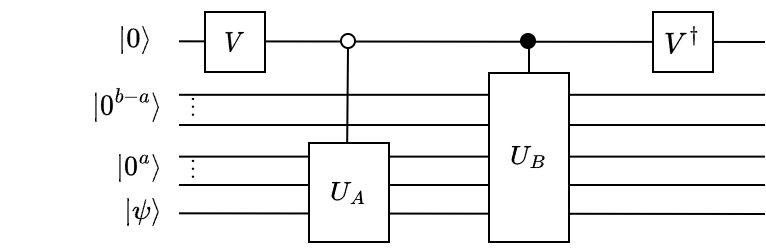}
	\renewcommand{\figurename}{Figure} 
	\caption{Circuit for addition of two block encoded matrices ($a\leq b$).}
	\label{fig:Addition01}
\end{figure} 
} 

\textbf{Theorem 4:}\label{Theorem:4} Multiplication of two block-encoded matrices. {If $U_A$ is an ($\alpha,a,\delta$)-block-encoding of an $s$-qubit operator $A$, and $U_B$ is a ($\beta,b,\epsilon$)-block-encoding of an $s$-qubit operator $B$, then $(I_b\otimes U_A)(I_a \otimes U_B)$ is an $(\alpha\beta, a+b, \alpha\epsilon + \beta\delta)$-block-encoding of \emph{AB} in Lamma 53 in reference \cite{Lin2022}.

The quantum circuit for block encoding of $AB$ can be constructed in Figure \ref{fig: Multiplication}.  
\begin{figure}[!h]
	\centering
	\includegraphics[scale=0.3]{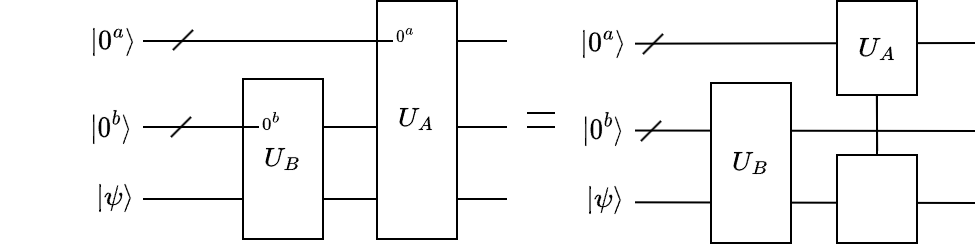}
	\renewcommand{\figurename}{Figure} 
	\caption{Circuit for product of two block encoded matrices.}
	\label{fig: Multiplication}
\end{figure} 
}

It should be noted that a marker is placed on the left side of the gate inspired by Qiskit to specify which ancilla qubits are affected by the quantum gate. For instance, we utilize $0^a$ as a marker to indicate that $U_A$ operates on $|0^a\rangle$ while no operation is performed on other qubits $|0^{b}\rangle$, and vice versa. Of course, both quantum operators act on state $|\psi\rangle$. This notation will be employed consistently throughout the rest of this article.

\subsection{Matrix Inversion by QSVT}

The HHL algorithm\cite{Harrow2009} is widely recognized as the most prominent quantum algorithm for matrix inversion. In the case of a sparse Hermitian matrix $A$, it enables us to efficiently solve the linear equation $Ax=b:$
\begin{align}
	|x\rangle=A^{-1}|b\rangle=\sum\limits_{i=1}^{n^2}\beta_i\frac{1}{\lambda_i}|u_i\rangle,
\end{align}
where $\lambda_i$ is the eigenvalue of the matrix $A$, $u_i$ is the corresponding eigenvector. However, we need to dilate the matrix into a Hermitian matrix $\tilde{A}=\begin{bmatrix}0 & A \\ A^\dagger & 0\\ \end{bmatrix}$ if $A$ is a general square matrix. Additionally, obtaining the explicit form of $A^{-1}$ poses significant challenges and hinders further calculation. In this paper, we consider addressing the problem of matrix inversion through the quantum singular value transformation (QSVT)\cite{Lin2022,Martyn2021}. 
	
Given an $s$-qubit square matirx $A\in\mathbb{C}^{n\times n}$ with the singular value decomposition $A=W \Sigma V^{\dagger}$, where $\Sigma$ is a diagonal matrix of the singular values. Referring to the setup of HHL, we assume that the singular values of $A$ satisfy $\sigma_i \in [\kappa^{-1},1]$, where $\kappa$ is condition number of $A$. As all the singular values are nonzero, we can obtain the inversion of $A$ as  $A^{-1}=V \Sigma^{-1} W^{\dagger}$, where $\Sigma^{-1}$ contains reciprocals of the singular values along the diagonal. Because $A^{\dagger}=V \Sigma W^{\dagger}$, then
\begin{align}
	A^{-1}=Vf(\Sigma)W^\dagger=f(A^\dagger)
\end{align}
where $f(x)=x^{-1}$ is an odd function. According to the QSVT with real polynomials, we can construct the circuit of $U_{p(A)}\in BE_{1,a+1}(p(A))$ when $f(x)$ is an odd function, as depicted in Figure \ref{fig: QSVT_d_odd}. Here, $p(x)$ represents an odd polynomial that approximates $f(x)$. Further details can be found in references \cite{Gilyen2019,Lin2022,Martyn2021}. 

\begin{figure}[!h]
	\centering
	\includegraphics[scale=0.4]{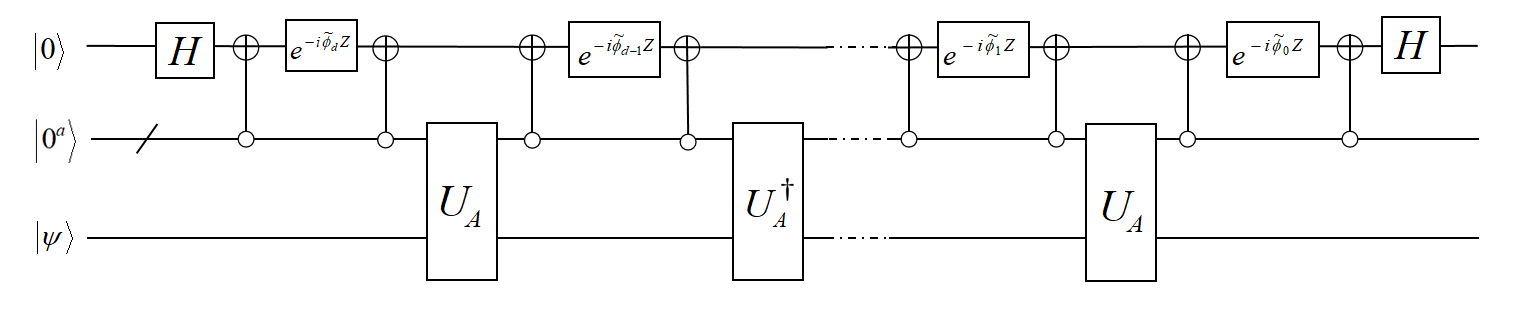}
	\renewcommand{\figurename}{Figure} 
	\caption{Circuit of $U_{P(A)}\in BE_{1,a+1}(P(A))$ when d is odd.}
	\label{fig: QSVT_d_odd}
\end{figure} 

Hence, it is imperative to identify an odd polynomial $p(x)$ that serves as an $\epsilon'$-approximation to the function $f(x)$ within the interval $[\frac{1}{\kappa},1]$. Notably, based on \cite{Lin2022}, we can obtain an odd polynomial $p(x)$ of degree $O(\kappa\log(1/\epsilon'))$ which effectively approximates $\frac{1}{\kappa\beta}\frac{1}{x}$. That is
\begin{align}
	|p(x)-\frac{1}{\kappa\beta}\frac{1}{x}|\leq \epsilon',\quad \forall x\in[-1,-\frac{1}{\kappa}]\cup[\frac{1}{\kappa},1].
	\label{equ:Matrix_Inversion_Poly}
\end{align}
where $p(x)$ has odd parity and $|p(x)|\le 1$ for $x\in[-1,-\frac{1}{\kappa}]\cup[\frac{1}{\kappa},1]$, $\beta$ is chosen arbitrarily.

The matrix inversion algorithm by QSVT is designed as follows\cite{Lin2022,Martyn2021}:

Given an $n \times n$ matrix $A$, where $\Vert A \Vert \leq 1$(if not, $A$ can be rescaled to meet this condition), and a precision parameter $\epsilon'$.

(1) Prepare the unitary block encoding of $A$ and $A^\dagger$, denoted by $U_A\in BE_{1,a}(A),U_{A^\dagger}\in BE_{1,a}(A^\dagger)$ respectively.

(2) Construct the circuit shown in Figure \ref{fig: QSVT_d_odd} to get $U_{p(A)}\in BE_{1,a+1}(p(A))$, where $p(x)$ represents the matrix inversion polynomial of degree $O(\kappa\log(1/\epsilon'))$. Note that the phase factor $\widetilde{\Phi}=(\widetilde{\phi_0},\widetilde{\phi_1},\cdots,\widetilde{\phi_d})$ can be determined from the polynomial $p(x)$.

(3) Finally, we have that $p(A)$ is an $\epsilon$-approximation to $\frac{1}{\kappa\beta}A^{-1}$.

\section{A quantum algorithm of the Kalman filter}

\subsection{The classical Kalman filter}

The Kalman filter is employed for state estimation in linear dynamical systems represented in the state space format\cite{Kim2018}. The process model specifies the transition of the state from time $k-1$ to time $k$ as follows:	
\begin{align}
	x_{k} &=Ax_{k-1}+Bu_{k-1}+w_{k-1}, \ k \geq 1.
	\label{equ:c_KF_state}
\end{align}
where $A$ is the state transition matrix applied to the previous state vector $x_{k-1}$, $B$ is the control-input matrix applied to the control vector $u_{k-1}$, and $w_{k-1}$ is the process noise vector which is zero-mean Gaussian with the covariance $Q$, i.e., it can be represented as $w_{k-1}\sim N(0,Q)$.  

The process model is paired with the observation model which characterizes the correlation between the state and the observation at time step $k$ as:
\begin{align}
	z_{k} &=Hx_{k}+v_{k}.
	\label{equ:c_KF_observation}
\end{align}
where $z_k$ is the observation vector, $H$ denotes the observation matrix, and $v_k$ represents the observation noise vector which is zero-mean Gaussian with the covariance $R$, i.e., $v_{k}\sim N(0,R)$. Let's assume that both $w_{k-1}$ and $v_{k}$ are statistically independent.

The Kalman filter algorithm for state estimation comprises two stages: prediction and update. Essentially, it employs the Kalman gain to refine the predicted state value towards its true value counterpart. The Kalman filter algorithm is summarized as follows.

1. The prior state estimation is calculated:
\begin{align}
	\hat{x}_{k}^-=A\hat{x}_{k-1}+Bu_{k-1}.
	\label{equ:c_KF01}
\end{align}

2. The prior error covariance is computed: 
\begin{align}
	P_{k}^-=AP_{k-1}A^T+Q.
	\label{equ:c_KF02}
\end{align}

3. The Kalman gain is determined: 
\begin{align}
	K_{k}=P_{k}^-H^T[HP_{k}^-H^T+R]^{-1}.
	\label{equ:c_KF03}
\end{align}

4. The state estimation is updated: 
\begin{align}
	\hat{x}_k=\hat{x}_k^-+K_{k}(z_k-H\hat{x}_k^-).
	\label{equ:c_KF04}
\end{align}

5. The error covariance is computed:
\begin{align}
	P_k=(I-K_kH)P_k^-.
	\label{equ:c_KF05}
\end{align}

Equations (\ref{equ:c_KF01})$\sim$(\ref{equ:c_KF02}) represent the discrete Kalman filter prediction equations, while equations (\ref{equ:c_KF03})$\sim$(\ref{equ:c_KF05}) correspond to the discrete Kalman filter update equations. This procedure ensures accurate estimation of the state variables $\hat{x}_k$. The whole process involves computation of these five formulas through matrix operations. Considering that block encoding of arbitrary matrices used in matrix operations can lead to exponential speed-ups\cite{Gilyen2019}, we propose a corresponding quantum algorithm for constructing the Kalman filter based on block encoding.

\subsection{The quantum algorithm}

The implementation procedure of the quantum algorithm for the Kalman filter involves evolving the encoded quantum state to the final required state through a series of unitary transformations, followed by obtaining the final result via quantum measurement. The process begins with block encoding and ends with measurement. 

Before we proceed, let me make two brief statements: 

(1)The matrices $A\in \mathbb{R}^{n \times n}, B\in \mathbb{R}^{n \times c},  H\in \mathbb{R}^{m \times n}, Q\in\mathbb{R}^{n \times n}, R \in \mathbb{R}^{m\times m}$ in equations (\ref{equ:c_KF_state})$\sim$(\ref{equ:c_KF_observation}) are known and constant, where $1\leq m,c,n$. At the initial time, the vectors $\hat{x}_0\in\mathbb{R}^n,u_0\in\mathbb{R}^c$ and the matrix $P_0\in\mathbb{R}^{n \times n}$ are given. Moreover, we assume that $u_0=u_1=\cdots=u_K$. It is considered that we can get $z_k\in\mathbb{R}^m$ from the sensor measurements at time $k$.

(2) To analyze the complexity of this quantum algorithm, we adopt certain criteria for block encoding: if an arbitrary $s$-qubit matrix $A\in \mathbb{C}^{n\times n}$ is stored in a quantum-accessible data structure, then we can always contruct a $(\Vert A\Vert_F,s,\epsilon)$-block-encoding of $A$ that can be implemented in time $O(polylog(n/\epsilon))$ refer to Theorem 2. For simplicity, let's assume that all these matrices have an error-free block encoding.

We present the complete procedure for implementing the Kalman filter using a quantum approach in Figure \ref{fig:algorithm_flow}. This includes addition, multiplication and inversion operations on the block-encoded matrices. The detailed process is given below. 
\begin{figure}[!h]
	\centering
	\includegraphics[scale=0.6]{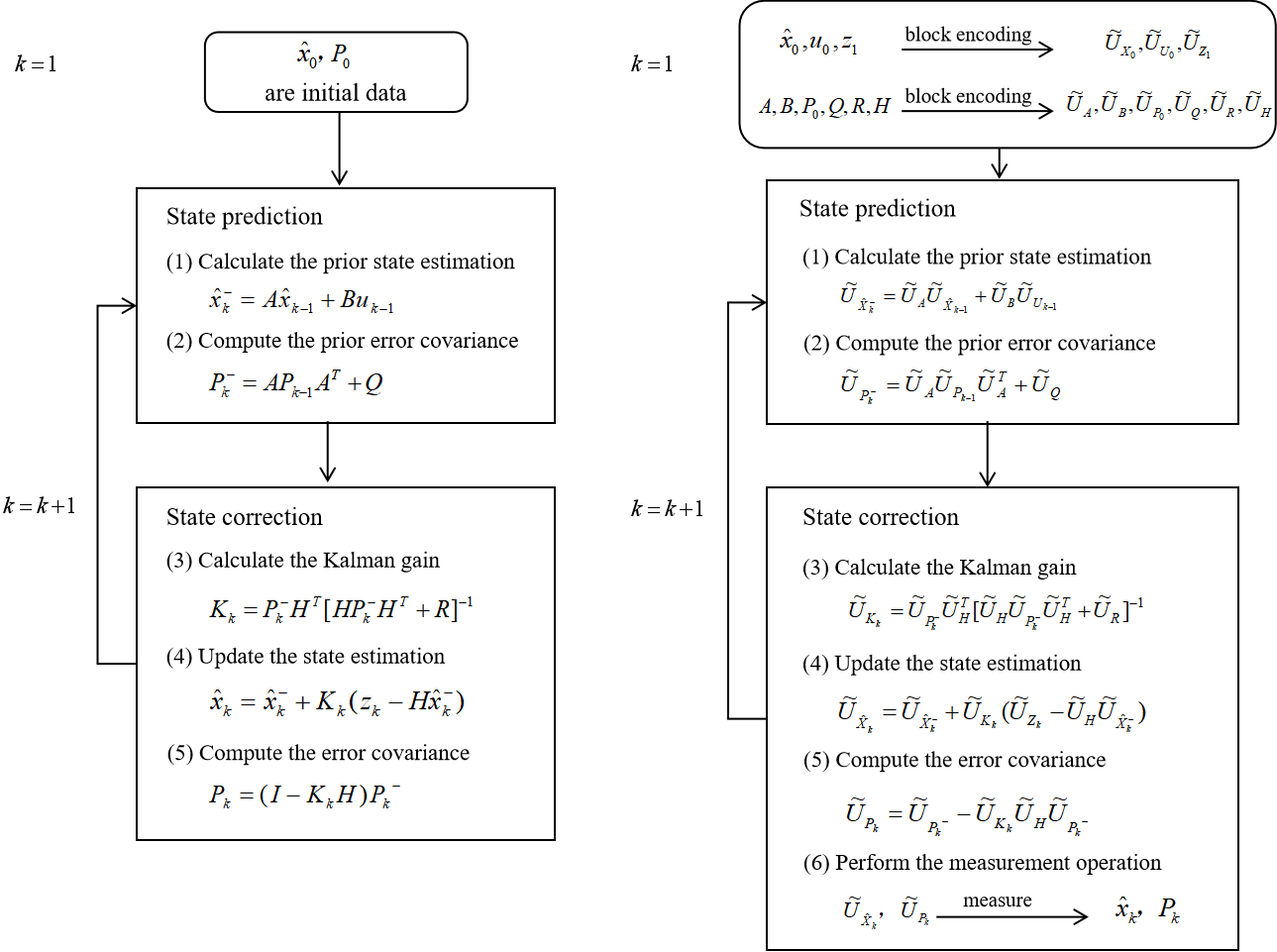}
	\renewcommand{\figurename}{Figure} 
	\caption{The Kalman filter (left: the classical algorithm, right: the quantum algorithm)}
	\label{fig:algorithm_flow}
\end{figure} 

\subsubsection{Block encoding of the known data}

Since $\hat{x}_0$,$u_0$,$z_1$ are vectors of dimensions $n$, $c$ and $m$ respectively, for simplicity of calculation, they need to be converted into $n\times n$ square matrices for block encoding. This can be achieved by filling in zeros in the remaining positions of matrices as follows:
\begin{align}
 	\hat{X}_0=\begin{bmatrix} 
 		\hat{x}_0^{(0)}  & 0 &\cdots &\cdots & 0 \\ 
 		\vdots & \vdots & \ddots & \ddots & \vdots \\
 		\vdots & \vdots & \ddots & \ddots & \vdots \\
 		\hat{x}_0^{(n-2)}  & 0 &\cdots &\cdots & 0 \\
 		\hat{x}_0^{(n-1)}  & 0 &\cdots &\cdots & 0 
 	\end{bmatrix}_{n\times n} 
 	U_0=\begin{bmatrix} 
 		u_0^{(0)}  & 0 &\cdots & \cdots & 0 \\
 		\vdots & \vdots & \ddots & \ddots &\vdots \\
 		u_0^{(c-1)}  & 0 &\cdots & \cdots & 0 \\
 		\vdots & \vdots & \ddots & \ddots &\vdots \\
 		0  & 0 &\cdots & \cdots & 0 
 	 \end{bmatrix}_{n\times n} 
  	Z_1=\begin{bmatrix}
  		z_0^{(0)}  & 0 &\cdots & \cdots & 0 \\
  		\vdots & \vdots & \ddots & \ddots &\vdots \\
  		z_0^{(m-1)}  & 0 &\cdots & \cdots & 0 \\
  		\vdots & \vdots & \ddots & \ddots &\vdots \\
  		0  & 0 &\cdots & \cdots & 0  
  	\end{bmatrix}_{n\times n}  
\end{align}
The block encoding of a matrix is denoted by $\widetilde{U}$ in order to differentiate it from the control matrix $U_{k-1}$.

Since $B,H$ are not necessarily square matrices and $R$ is not an $n$-order square matrix, it is necessary to convert them into $n\times n$ matrices for block encoding by padding the remaining positions with zeros.

We can perform block encoding for all known matrices using Theorem 2 as shown in Table \ref{tab:table1}.
\begin{table}[!t]
	\renewcommand\arraystretch{2}
	\caption{The block encoding of the known data\label{tab:table1}}
	\centering
	\footnotesize
	\begin{tabular}{ccc}
		\hline
		Original matrix &  Dimension & Block encoding  \\
		\hline
		$A$ & $n\times n$ & $\widetilde{U}_A \in BE_{\Vert A\Vert_F,s}(A)$\\
		\hline
		$B$ & $n\times c$  & $\widetilde{U}_B \in BE_{\Vert B\Vert_F,s}(B)$\\
		\hline
		$\hat{X}_0$ & $n\times n$  & $\widetilde{U}_{\hat{X}_0}\in BE_{\Vert\hat{X}_0 \Vert_F,s}(\hat{X}_0)$\\
		\hline
		$U_0$ & $n\times n$  & $\widetilde{U}_{U_0}\in BE_{\Vert U_{k-1}\Vert_F,s}(U_0)$\\
		\hline
		$P_0$ & $n\times n$  & $\widetilde{U}_{P_0}\in BE_{\Vert P_{k-1}\Vert_F,s}(P_0)$\\
		\hline
		$Q$ & $n\times n$  & $\widetilde{U}_Q\in BE_{\Vert Q\Vert_F,s}(Q)$\\
		\hline
		$R$ & $n\times n$  & $\widetilde{U}_R\in BE_{\Vert R\Vert_F,s}(R)$\\
		\hline
		$H$ & $m\times n$  & $\widetilde{U}_H\in BE_{\Vert H\Vert_F,s}(H)$\\
		\hline
		$Z_1$ & $n\times n$  & $\widetilde{U}_{Z_1}\in BE_{\Vert Z_1\Vert_F,s}(Z_1)$\\
		\hline
	\end{tabular}
\end{table}

\subsubsection{The computational process of the quantum algorithm of the Kalman filter}

We summarize the procedure of computing the Kalman filter based on the block encoding framework in Algorithm 1. The results of each procedure are shown in Table \ref{tab:table2}. According to Theorems 3 and 4, we can successfully perform the addition and multiplication operations on unitary matrices. More importantly, it is necessary to solve for the inverse of matrix $\widetilde{U}_{A_{temp}}$, which can be achieved through the QSVT in procedure 5. To obtain $A_{temp}$, we first need to measure $\widetilde{U}_{A_{temp}}$ to obtain $A_{temp}$ as it provides essential a priori knowledge such as the condition number $\kappa$ and singular values $\sigma$. This information enables us to calculate suitable phase factors for approximating $1/x$. Immediately according to equation (\ref{equ:Matrix_Inversion_Poly}), the $\epsilon'$ polynomial approximating $f(x)=\frac{1}{\kappa\beta}\frac{1}{x}$ can be expressed as $p(x)$. Then unitary operator $\widetilde{U}_{A_{temp}^{-1}}$ can be constructed. In the end, if the loop continues, we need to reconstruct the block encoding of $\hat{X}_k$, $P_k$ and $Z_k$ for the next calculation. Using the $\widetilde{U}_{\hat{X}_k}, \widetilde{U}_{P_k}$ directly leads to more qubits being used later. Moreover, the more ancilla qubits there are in construction of $U_A$, the larger normalization factor $\alpha$ is, it's more difficult to obtain the exact matrix $A$ by measurement.

The quantum circuit for a loop is illustrated in Figures \ref{fig:X_k} and \ref{fig:P_k}. Let's provide a brief explanation of the symbols depicted in these figures. As shown in Figure \ref{fig:X_k}, $\widetilde{U}$ represents the block encoding of the matrix, with its subscript denoting the given matrix. According to Theorem 3, operator $V$ is determined by the coefficients of two adder numbers. For instance, during the implementation process of $\widetilde{U}_{\hat{X}_k}$, both matrices $\widetilde{U}_{\hat{X}^-_k}$ and $\widetilde{U}_{K_k}(\widetilde{U}_{Z_k}-\widetilde{U}_H\widetilde{U}_{\hat{X}_k^-})$ that need to be summed are multiplied by their respective coefficients, which serve as scaling factors for the matrix. Operator $V_{\hat{X}_k}$ can be computed using these two coefficients. Similarly, the operator $V_{temp}, V_{P_k^-}, V_{\hat{X}^-_k}$ can be obtained following a similar procedure where we assign $temp=Z_k-H{\hat {X}^-_k}$.

Here, it is particularly noted that $\widetilde{U}_{A_{temp}^{-1}}$ is obtained by QSVT. Please refer to Figures \ref{fig:A_temp_inversion} and \ref{fig:A_temp} for the circuits of calculating $\widetilde{U}_{A_{temp}^{-1}}$ and $\widetilde{U}_{A_{temp}}$, respectively.

\begin{algorithm}[!h]  
\caption{The process of solving Kalman filter based on the block encoding framework.}\label{alg:alg1}  
\begin{algorithmic} 
        \State \hspace{0.3cm}\textbf{Input:} Access to all the given matrices $A,B,\hat{X}_0,U_0,P_0,Z_1,Q,R,H$.
		\State \hspace{0.3cm}\textbf{Output:} $\hat{X}_K, P_K$.
		\State \hspace{0.3cm}\textbf{Runtime:} $O(poly\log(n/\epsilon)\kappa\log(1/\epsilon')) $, the time complexity is analyzed in detail in Section 4.2.
		\State
		\State \hspace{0.3cm}\textbf{Procedure:}
		\State 1 \hspace{0.1cm}Prepare the unitary block encoding of $A,B,\hat{X}_0,U_0,P_0,Z_1,Q,R,H$.
		\State 2 \hspace{0.1cm}\textbf{for} $k=1$ to $K$ \textbf{do}:
		\State
		\State 3 \hspace{0.5cm}Calculate the block encoding of the prior state $\widetilde{U}_{\hat{X}_{k}^-}=\widetilde{U}_A \widetilde{U}_{\hat{X}_{k-1}} + \widetilde{U}_B \widetilde{U}_{U_{k-1}}$.
		\State 4 \hspace{0.5cm}Compute the block encoding of the prior error covariance  $\widetilde{U}_{P_k^-}=\widetilde{U}_A\widetilde{U}_{P_{k-1}}\widetilde{U}_A^T+\widetilde{U}_Q$.
		\State 5 \hspace{0.5cm}Calculate the block encoding of the Kalman gain $\widetilde{U}_{K_k}=\widetilde{U}_{P_k^-}\widetilde{U}_H^T[\widetilde{U}_H\widetilde{U}_{P_k^-}\widetilde{U}_H^T+\widetilde{U}_R]^{-1}$, where \\     \hspace{0.8cm}$\widetilde{U}_{A_{temp}}=\widetilde{U}_H\widetilde{U}_{P_k^-}\widetilde{U}_H^T+\widetilde{U}_R$. Measure $\widetilde{U}_{A_{temp}}$ to get $A_{temp}$, which is subsequently subjected  to\\
		\hspace{0.7cm} block encoding again after knowing the condition number and the singular values.
		\State 6 \hspace{0.5cm}Update the block encoding of the posterior state estimation  $\widetilde{U}_{\hat{X}_k}=\widetilde{U}_{\hat{X}_k^-}+\widetilde{U}_{K_k}(\widetilde{U}_{Z_k}-\widetilde{U}_H\widetilde{U}_{\hat{X}_k^-})$.
		\State 7 \hspace{0.5cm}Update the block encoding of the error covariance $\widetilde{U}_{P_k}=(\widetilde{U}_I-\widetilde{U}_{K_k}\widetilde{U}_H)\widetilde{U}_{P_k^-}$.
		\State 8 \hspace{0.5cm}Conducting measurements on $\widetilde{U}_{\hat{X}_k}, \widetilde{U}_{P_k}$ to obtain $\hat{X}_k, P_k$.
		\State 9 \hspace{0.5cm}Block encode $\hat{X}_k, P_k, Z_k$
		\State
		\State 10 \hspace{0.1cm} If the loop ends, the final results are $\hat{X}_K, P_K$.
	\end{algorithmic}
	\label{alg1}
\end{algorithm}

\begin{table}[!h]
	\renewcommand\arraystretch{2}
	\caption{The result of the intermediate variables\label{tab:table2}}
	\centering
	\footnotesize
	\begin{tabular}{ccc}
		\hline
		Procedure & Result &  Normalization factor \\
		\hline
		\multirow{3}{*}{3}  &	$\widetilde{U}_{A}\widetilde{U}_{\hat{X}_{k-1}}\in BE_{\alpha_{31},2s}(A\hat{X}_{k-1})$ &  $\alpha_{31}=\Vert A\Vert_F\Vert\hat{X}_{k-1}\Vert_F$ \\
		
		&$\widetilde{U}_{B}\widetilde{U}_{U_{k-1}}\in BE_{\alpha_{32},2s}(BU_{k-1})$ &  $\alpha_{32}=\Vert B\Vert_F\Vert\hat{X}_{k-1} \Vert_F$ \\
		
		&$\widetilde{U}_{\hat{X}_{k}^-}\in BE_{\alpha_{\hat{X}_k^-},2s+1}(\hat{X}_k^-)$ &  $\alpha_{\hat{X}_k^-}=\alpha_{31}+\alpha_{32}$\\
		\hline
		\multirow{2}{*}{4} & $\widetilde{U}_A\widetilde{U}_{P_{k-1}}\widetilde{U}_{A^T}\in BE_{\alpha_{41},3s}(A P_{k-1}A^T)$  & $\alpha_{21}={\Vert A\Vert}_F^2 \Vert P_{k-1}\Vert_F$ \\
		
		& $\widetilde{U}_{P_k^-}\in BE_{\alpha_{P_k^-},3s+1}(P_k^-)$ & $\alpha_{P_k^-}=\alpha_{41}+\Vert Q\Vert_F$\\
		\hline
		\multirow{5}{*}{5} & $\widetilde{U}_{P_k^-}\widetilde{U}_{H^T}\in BE_{\alpha_{51},4s+1}(P_k^-H^T)$  & $\alpha_{51}=\alpha_{P_k^-}\Vert H\Vert_F$ \\
		
		& $\widetilde{U}_H\widetilde{U}_{P_k^-}\widetilde{U}_{H^T}\in BE_{\alpha_{52},5s+1}(H P_k^-H^T)$ & $\alpha_{52}=\alpha_{51}\Vert H\Vert_F$\\
		
		& $\widetilde{U}_{A_{temp}}\in BE_{\alpha_{53},5s+2}(A_{temp})\xrightarrow[]{block encoding}BE_{\alpha_{53}',s}(A_{temp})$ & $\alpha_{53}=\alpha_{52}+\Vert R\Vert_F, \alpha_{53}'=\gamma\alpha_{53}$ \\
		
		& $\widetilde{U}_{A^{-1}_{temp}}\in BE_{\alpha_{54},s+1}(A_{temp}^{-1},\epsilon')$ & $\alpha_{54}=\kappa\beta/\alpha_{53}'$ \\
		
		& $\widetilde{U}_{K_k}\in BE_{\alpha_{K_k},5s+2}(K_k)$ & $\alpha_{K_k}=\kappa\beta\alpha_{51}/\alpha_{53}'$ \\
		\hline
		\multirow{4}{*}{6} & $-\widetilde{U}_H\widetilde{U}_{\hat{X}_k^-}\in BE_{\alpha_{61},3s+1}(-H\hat{X}_k^-)$  & $\alpha_{61}=\alpha_{\hat{X}_k^-}\Vert H\Vert_F$ \\
		
		& $\widetilde{U}_{Z_k}-\widetilde{U}_H\widetilde{U}_{\hat{X}_k^-}\in BE_{\alpha_{62},3s+2}(Z_k-H\hat{X}_k^-)$ & $\alpha_{62}=\alpha_{61}+\Vert Z_k\Vert_F$ \\
		
		& $\widetilde{U}_{K_k}(\widetilde{U}_{Z_k}-\widetilde{U}_H\widetilde{U}_{\hat{X}_k^-})\in BE_{\alpha_{63},8s+4}(K_k(Z_k-H\hat{X}_k^-))$ & $\alpha_{63}=\alpha_{62}\alpha_{K_k}$ \\
		
		& $\widetilde{U}_{\hat{X}_k}\in BE_{\alpha_{\hat{X}_k},8s+5}(\hat{X}_k)$ & $\alpha_{\hat{X}_k}=\alpha_{\hat{X}_k^-}+\alpha_{63}$ \\
		\hline
		\multirow{2}{*}{7} & $\widetilde{U}_{K_k}\widetilde{U}_H\widetilde{U}_{P_k^-}\in BE_{\alpha_{71},9s+3}(K_kHP_k^-)$  & $\alpha_{71}=\alpha_{K_k}\alpha_{P_k^-}\Vert H\Vert_F$ \\
		
		& $\widetilde{U}_{P_k} \in BE_{\alpha_{P_k},9s+4}(P_k)$ &  $\alpha_{P_k}=\alpha_{71}+\alpha_{P_k^-}$ \\
		\hline
		\multirow{2}{*}{8} & $\hat{X}_k=\alpha_{\hat{X}_k}(\langle 0^{8s+5}|\otimes I_s)\widetilde{U}_{\hat{X}_k}(|0^{8s+5}\rangle \otimes I_s)$  & $\alpha_{\hat{X}_k}=\alpha_{\hat{X}_k^-}+\alpha_{63}$ \\
		
		& $P_k=\alpha_{P_k}(\langle 0^{9s+4}|\otimes I_s)\widetilde{U}_{P_k}(|0^{9s+4}\rangle \otimes I_s)$ &  $\alpha_{P_k}=\alpha_{71}+\alpha_{P_k^-}$ \\
		\hline
	\end{tabular}
\end{table}

\begin{figure}[!h]
	\centering
	\includegraphics[scale=0.5]{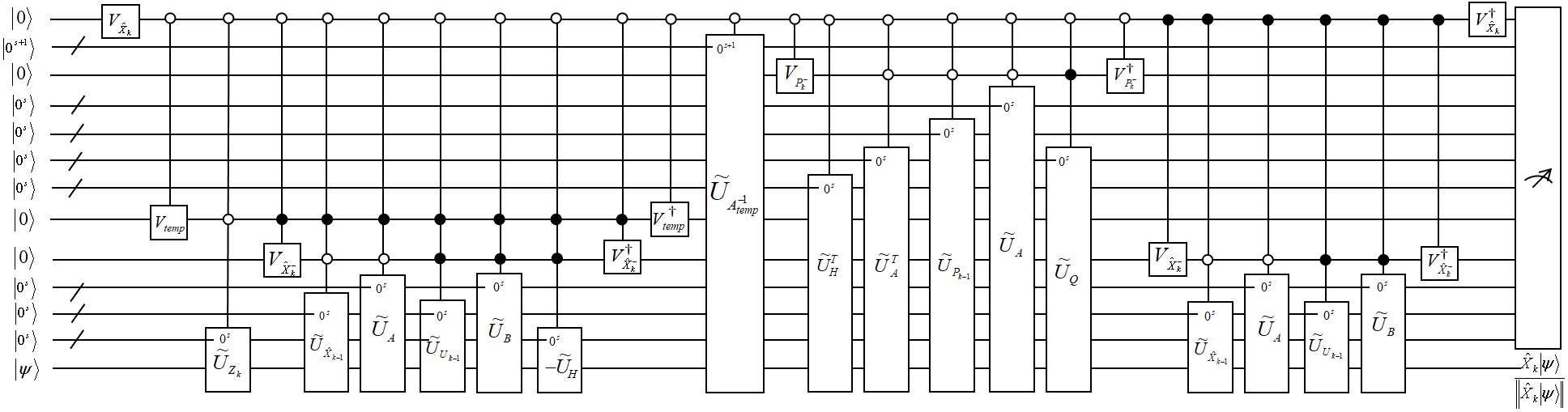}
	\renewcommand{\figurename}{Figure} 
	\caption{The quantum circuit of calculating $\hat{X}_K$}
	\label{fig:X_k}
\end{figure} 

\begin{figure}[!h]
	\centering
	\includegraphics[scale=0.55]{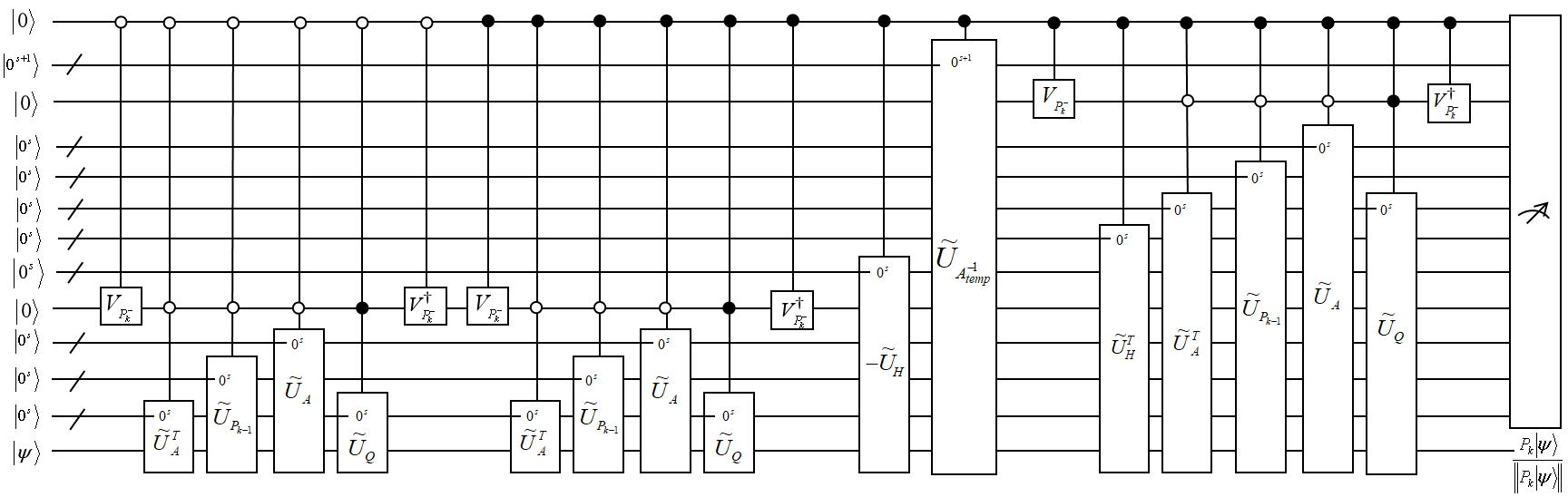}
	\renewcommand{\figurename}{Figure} 
	\caption{The quantum circuit of updating $P_K$}
	\label{fig:P_k}
\end{figure} 

\begin{figure}[!h]
	\centering
	\includegraphics[scale=0.8]{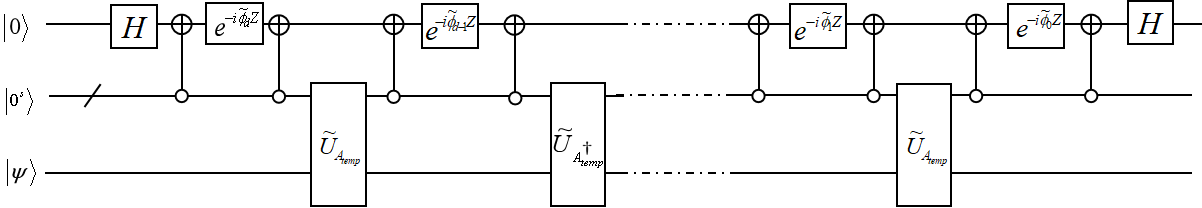}
	\renewcommand{\figurename}{Figure} 
	\caption{The quantum circuit of calculating $\widetilde{U}_{A_{temp}^{-1}}$}
	\label{fig:A_temp_inversion}
\end{figure}

\begin{figure}[!h]
	\centering
	\includegraphics[scale=0.4]{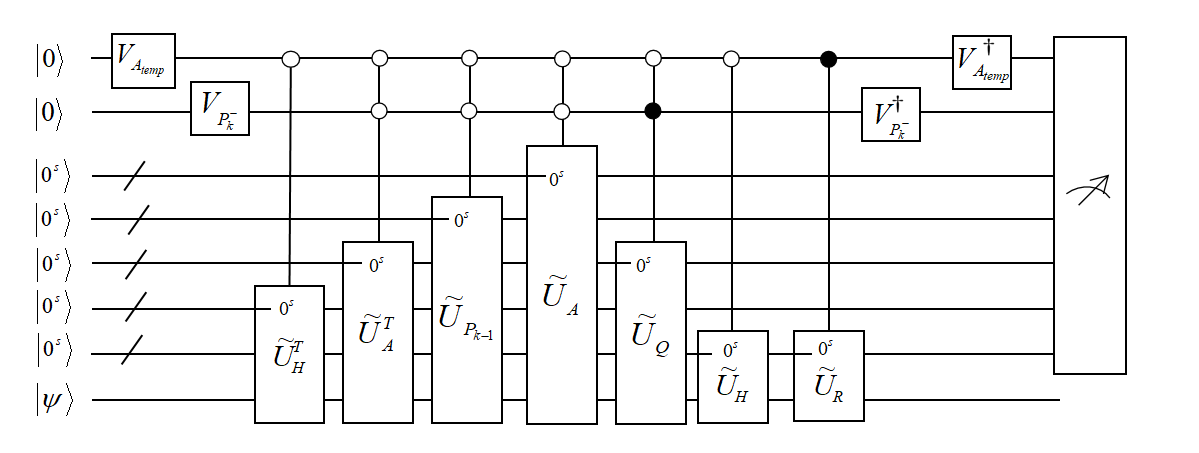}
	\renewcommand{\figurename}{Figure} 
	\caption{The quantum circuit of calculating $\widetilde{U}_{A_{temp}}$}
	\label{fig:A_temp}
\end{figure}

\section{Algorithm complexity analysis}
In this section, we conduct an analysis of the time and space complexities involved.

\subsection{Time complexity theory}

Now, let's analyse the time complexity of this algorithm, which comprises four parts: block encoding, addition of unitary matrices, multiplication of unitary matrices, and matrix inversion. 

(1) The time complexity of block encoding. 

According to Theorem 1, if $A$ is an $s$-qubit square matrix, where $A \in \mathbb{C}^{n\times n}, s=\lceil\log n\rceil$, the quantum-accessible data structure stores the matrix $A$ in times $O(\log^2 n)$\cite{Chakraborty2018}, Then $A$ can be efficiently block encoded to an $\epsilon$-approximate unitary $U_A$, an $(\Vert A\Vert_F,s,\epsilon)$-block encoding of $A$, with a time complexity of $O(polylog(n/\epsilon))$ in the quantum-accessible data structure. So the quantum-accessible data structure enables the transformation of $A$ to $U_A$ in a time complexity of $O(polylog(n/\epsilon))$, denoted as $T_{U_A}$.

(2) The time complexity of matrix addition.

If $U_A$ is an ($\alpha,a,\epsilon_a$)-block-encoding of an $s$-qubit operator $A$ that can be implemented in time $T_{U_A}$, and if $U_B$ is a ($\beta,b,\epsilon_b$)-block-encoding of an $s$-qubit operator $B$ that can be implemented in time $T_{U_B}$, then according to Theorem 3, we can perform the $(\alpha+\beta,\max\{a,b\}+1,\epsilon_a+\epsilon_b)$-block-encoding of the sum of operators A and B in time complexity $O(T_{U_A}+T_{U_B})$.

(3) The time complexity of matrix multiplication.

With $U_A$ and $U_B$ give above, according to Theorem 4, we can perform the $(\alpha+\beta,\max\{a,b\}+1,\epsilon_a+\epsilon_b)$-block-encoding of $AB$ in time complexity O($T_{U_A}+T_{U_B}$).

(4) The time complexity of matrix inversion.

To solve $A^{-1}$, we introduce the matrix polynomial $p(A)$ of degree $O(\kappa\log(1/\epsilon'))$, which serves as an $\epsilon'$-approximation to $(\frac{1}{\kappa\beta})\frac{1}{A}$. Here, $\kappa$ represents the condition number and $\epsilon'$ denotes the desired precision. The degree $d=O(\kappa\log(1/\epsilon'))$ implies that QSVT requires implementing $O(d)$ instances of $U_A,U_A^\dagger$ and $(d+1)$ interleaved rotations, as depicted in Figure \ref{fig: QSVT_d_odd}. Consequently, the matrix inversion process for obtaining $A^{-1}$ through QSVT can be accomplished within a time complexity of $O(T_{U_A}d)$.

\subsection{Algorithm complexity analysis}

\subsubsection{Time complexity analysis}

We analyze the time complexity of each procedure and compare it with the classical algorithm. 

(1) Block encoding. Before performing block encoding, all the matrices need to have been stored in quantum-accessible data structure, then the block encoding of $A,B,X_{k-1},\\U_{k-1},P_{k-1},Q,R,Z_k$ in quantum-accessible data structure can be implemented, the time cost of each matrix are all the same which can be denoted $T_U$, so the block encoding process is finished in time $O(T_U)$, where $T_U=poly\log(n/\epsilon)$. 

(2) For $\widetilde{U}_{\hat{X}_k^-}$, the total computation time is 4$T_U$, thus enabling the implementation of $\widetilde{U}_{\hat{X}_k^-}$ in $O(T_U)$ time. Similarly, we can also accomplish $\widetilde{U}_{P_k^-}$ within $O(T_U)$ time.

(3) Calculating $\widetilde{U}_{K_k}$ involves addition, multiplication and invertion operations. The total time required for matrix inversion $\widetilde{U}_{A_{temp}^{-1}}$ is $7T_Ud$, where $d=O(\kappa\log(1/\epsilon'))$, with $\epsilon'$ representing the desired precision of the solution and $\kappa$ denoting the condition number ($\kappa >1$). Consequently, the time complexity of this process can be expressed as $O(T_Ud)=O(poly\log(n/\epsilon)\kappa\log(1/\epsilon'))$.

(4) The calculation step for $\widetilde{U}_{\hat{X}_k^-}$ involves the operation of $\widetilde{U}_{K_k}$, enabling the simulation of $\widetilde{U}_{\hat{X}_k^-}$ in a time complexity of $O(poly\log(n/\epsilon)\kappa\log(1/\epsilon'))$, as well as $\widetilde{U}_{P_k}$.

In conclusion, the overall time complexity of the algorithm is $O(poly\log(n/\epsilon)\kappa\log(1/\epsilon'))$.

Next, we present a complexity comparison between this algorithm and its classical counterpart. In the process of solving the Kalman filter problem, various matrix operations need to be handled. With brute force computation, the holistic time complexity of the classical algorithm is $O(n^3)$. Specifically, naive matrix multiplication for square matrices with dimensions $n\times n$ has a complexity $O(n^3)$, while Strassen's algorithm achieves a lower complexity of  $O(n^{2.807})$\cite{Golub2012}. Additionally, there are several approaches reducing the computational complexity of the Kalman filter. Parallel computing can effectively decrease computation time by parallelizing matrix operations within the linear algebra library that handles these operations\cite{Raitoharju2019}. This leads to a reduction in time complexity from $O(n)$ to $O(\log n)$, in terms of the number of time steps taken, rather than reducing it per single step alone. A method known as Fast Algorithm for Kalman Filtering\cite{Kailath2000} can further reduce the computational complexity of $P_k$ from $O(n^3)$ to $O(n^2\alpha)$. From Table \ref{tab1}, we observe that the quantum algorithm exhibits a time complexity of approximately 	$O(poly\log(n/\epsilon)\kappa\log(1/\epsilon'))$, which still demonstrates exponential acceleration compared with its classical counterpart.

\begin{table}[!t]
\footnotesize
\centering
\renewcommand\arraystretch{2}
\caption{The time complexity of classical algorithm and quantum algorithm}
\label{tab1}
\tabcolsep 15pt 
\begin{center}
	\begin{tabular*}{\textwidth}{ccc}
		\toprule
		Procedure & Classical algorithm & Quantum algorithm \\
		\hline
		Block encoding & 0 & $O(polylog(n/\epsilon))$ \\
		$\hat{x}_{k}^-=A\hat{x}_{k-1}+Bu_{k-1}$ & $O(n^2)$ &$O(polylog(n/\epsilon))$ \\
		$P_{k}^-=AP_{k-1}A^T+Q$ & $O(n^3)$ & $O(polylog(n/\epsilon))$ \\
		$K_{k}=P_{k}^-H^T[HP_{k}^-H^T+R]^{-1}$ & $O(mn^2)$ & $O(poly\log(n/\epsilon)\kappa\log(1/\epsilon'))$ \\
		$\hat{x}_k=\hat{x}_k^-+K_{k}(z_k-H\hat{x}_k^-)$ & $O(mn)$ & $O(poly\log(n/\epsilon)\kappa\log(1/\epsilon'))$ \\
		$P_k=(I-K_kH)P_k^-$ & $O(n^3)$ & $O(poly\log(n/\epsilon)\kappa\log(1/\epsilon'))$ \\
		Sum & $O(n^3)$ & $O(poly\log(n/\epsilon)\kappa\log(1/\epsilon'))$ \\
		\bottomrule
	\end{tabular*}
\end{center}
\end{table}

\subsubsection{Space complexity analysis}

As is depicted in Figure \ref{fig:X_k} and Figure \ref{fig:P_k}, in addition to the qubits required for constructing $|\psi\rangle$, denoted by $s$, the number of the ancilla qubits required for block encoding as well as for implementing matrix addition and multiplication operations are $8s+5$ and $9s+4$. Additionally, 2-qubit or 3-qubit controlled-U (CU) gates can be found, where U comprises a combination of CNOT or Toffoli gates. The CU gate may involve up to 5-qubit control gates, which could require an additional set of 4 ancilla qubits for decomposition into elementary gates at least. It can be seen that the total number of qubits required is linearly dependent on $s$.

For instance, when considering an $8$-dimensional column vector $x_0$ which can be represented by 3 qubits, the corresponding ancilla qubits required for block encoding and matrix operations amount to 29 and 31 resepctively ($s$=3), excluding the ancilla qubits needed for some multi-qubit control gates decomposition. We propose two approaches to reduce the number of ancilla qubits used. 

The first approach begins with the construction of block encoding. As is well-known, for any $s$-qubit matrix $A$ satisfying $\Vert A\Vert_F\leq 1$, its singular value decomposition can be denoted as $W\Sigma V^\dagger$, where all the singular values in the diagonal matrix $\Sigma$ belong to the interval $[0,1]$. Subsequently, an $(s+1)$-qubit unitary matrix can be constructed. 

\begin{align}
	U_A =\begin{bmatrix} 
		W & 0 \\ 
		0 & I_s 
	\end{bmatrix}
	\begin{bmatrix}
		\Sigma & \sqrt{I_s-\Sigma^2}\\ 
		\sqrt{I_s-\Sigma^2}  & -\Sigma
	\end{bmatrix}
	\begin{bmatrix}
		V^\dagger & 0 \\ 
		0  & I_s
	\end{bmatrix}= 
	\begin{bmatrix}
		A & W\sqrt{I_s-\Sigma^2} \\ 
		\sqrt{I_s-\Sigma^2}V^\dagger  & -\Sigma
	\end{bmatrix}.
\end{align}
The matrix $A$ can be block encoded by a single ancilla qubit, denoted as $U_A$, which is a $(1,1)$-block-encoding\cite{Dongrandom2021}. However, it should be noted that achieving this in practice remains challenging. If such encoding is feasible, the number of required ancilla qubits would be reduced. Nevertheless, currently there is no known general method to construct an arbitrary metrix using only one ancilla qubit. Further research and development are still necessary in this area.

The second approach involves measuring each precedure of the algorithm once and subsequently re-block-encoding the resulting matrix. As an example, consider $\widetilde{U}_{\hat{X}_k^-}$ which requires $(2s+1)$ ancilla qubits. If we measure $\widetilde{U}_{\hat{X}_k^-}$ to obtain $\hat{X}_k^-$ and perform block encoding again, only $s$ ancilla qubits are needed. This choice proves advantageous, however, finding the block encoding for an arbitrary matrix is challenging without a general construction method available in practice. Moreover, errors can occur during quantum computing measurements, with continuous measurement leading to error accumulation. Hence, it is crucial to minimize intermediate measurement steps.

\section{An illustrative example}
The proposed quantum algorithm is demonstrated through a simple illustrative example. We have implemented the entire algorithm using Qiskit, an open-source software development kit designed for quantum computing\cite{Qiskit2023}. To ensure simplicity, we consider a state space representation as depicted below:
\begin{eqnarray}
	x_k =\begin{bmatrix} 
		1 & -1 \\ 
		1 & 1 
	\end{bmatrix}x_{k-1}+
	\begin{bmatrix}
		1  \\ 
		1  
	\end{bmatrix}u_{k-1}+
	w_{k-1} 
\end{eqnarray}
\begin{eqnarray}
	z_k =\begin{bmatrix}
		2 & 0 \\
		0 & 1
	\end{bmatrix}x_k+v_k
\end{eqnarray}
Initialize the initial state estimation as
$\hat{x}_0 =\begin{bmatrix} 2  \\ 1  \end{bmatrix}$, the initial prior covariance matrix as
$P_0=\begin{bmatrix}
	1 & 0 \\ 
	0 & 1 
\end{bmatrix}$. The process noise is modeled as 
$w\sim N(0,\begin{bmatrix}1 & 0 \\ 0 & 1 \end{bmatrix})$, while the observation noise is denoted by $v\sim N(0,1)$. Therefore, we have matrices 	
$A=\begin{bmatrix}1 & -1 \\ 1 & 1 \end{bmatrix}, B=\begin{bmatrix}1 \\ 1 \end{bmatrix}, Q=\begin{bmatrix}1 & 0 \\ 0 & 1 \end{bmatrix}, R=\begin{bmatrix}1 & 0 \\ 0 & 1 \end{bmatrix}, H=\begin{bmatrix}2 & 0 \\ 0 & 1 \end{bmatrix}$. The initial condition is set as $u_0 = 1$, and the vector $z_1$ is defined as $[1, 1]^T$ for the sake of simplicity.

The quantum algorithm of the Kalman filter is outlined as follows:

(1) Construct the block encoding of $A, B, P_0, Q, R, H, \hat{x}_0, u_0, z_1$, where $s$=1. These variables are represented as follows:
\begin{equation}
	\begin{aligned}
		\widetilde{U}_{A}=
		\begin{bmatrix}
			\frac{1}{2} & -\frac{1}{2} & * & *\\ 
			\frac{1}{2} &  \frac{1}{2} & *  & *\\
			* & * & * & * \\
			* & * & * & *
		\end{bmatrix}_{4\times 4}
		\widetilde{U}_{B}=
		\begin{bmatrix}
			\frac{1}{\sqrt{2}} & 0 & * & *\\ 
			\frac{1}{\sqrt{2}} & 0 & *  & *\\
			* & * & * & * \\
			* & * & * & *
		\end{bmatrix}_{4\times 4}
		\widetilde{U}_{P_0}=\widetilde{U}_{Q}=\widetilde{U}_{R}=
		\begin{bmatrix}
			\frac{1}{\sqrt{2}} & 0 & * & *\\ 
			0 & \frac{1}{\sqrt{2}} & *  & *\\
			* & * & * & * \\
			* & * & * & *
		\end{bmatrix}_{4\times 4} \\
		\widetilde{U}_{U_0}=
		\begin{bmatrix}
			1 & 0 & * & *\\ 
			0 & 0 & *  & *\\
			* & * & * & * \\
			* & * & * & *
		\end{bmatrix}_{4\times 4} 
		\widetilde{U}_{H}=
		\begin{bmatrix}
			\frac{2}{\sqrt{5}} & 0 & 0 & *\\ 
			0 & \frac{1}{\sqrt{5}} & *  & *\\
			* & * & * & * \\
			* & * & * & *
		\end{bmatrix}_{4\times 4}
		\widetilde{U}_{\hat{X}_0}=
		\begin{bmatrix}
			\frac{2}{\sqrt{5}} & 0 & * & *\\ 
			\frac{1}{\sqrt{5}} & 0 & * & *\\
			* & * & * &  * \\
			* & * & * &  *
		\end{bmatrix}_{4\times 4}
		\widetilde{U}_{Z_1}=
		\begin{bmatrix}
			\frac{1}{\sqrt{2}} & 0 & * & *\\ 
			\frac{1}{\sqrt{2}} & 0 & * & *\\
			* & * & * & * \\
			* & * & * & *
		\end{bmatrix}_{4\times 4} 
		\label{equ:block_encoding01}
	\end{aligned}
\end{equation}
We must ensure that all matrices are unitary matrices, thus requiring the multiplication of matrix elements by a normalization factor. Additionally, we assume an exact block encoding with $\epsilon=0$. Therefore, the block encoding of $A$ is denoted as $BE_{2,1}(A)$, and similar treatment is applied to other cases.

(2) Let's compute $\widetilde{U}_{A\hat{X_0}}$ and $\widetilde{U}_{BU_0}$ separately, then add them together to get $U_{\hat{X}_1^-}$. The detailed calculation procedure is shown as follows:
\begin{align}
	\widetilde{U}_{A\hat{X}_0}=
	\begin{bmatrix}
		\frac{1}{2\sqrt{5}} & 0 & * & \cdots & *\\  
		\frac{3}{2\sqrt{5}} & 0 & * & \cdots & *\\ 
		* & * & * & \cdots & * \\
		\vdots & \vdots & \vdots &  \ddots &  \vdots\\
		* & * & * & \cdots & *
	\end{bmatrix}\in BE_{2\sqrt{5},2}(A\hat{X}_0),
	\widetilde{U}_{BU_0}=
	\begin{bmatrix}
		\frac{1}{\sqrt{2}} & 0 & * & \cdots & *\\ 
		\frac{1}{\sqrt{2}} & 0 & * & \cdots & *\\ 
		* & * & * & \cdots & * \\
		\vdots & \vdots & \vdots &  \ddots &  \vdots\\
		* & * & * & \cdots & *
	\end{bmatrix}\in BE_{\sqrt{2},2}(BU_0).
\end{align}

Let $T=2\sqrt{5}\widetilde{U}_{A\hat{X}_0}+\sqrt{2}\widetilde{U}_{BU_0}$ to be LCU, so the block encoding of $\hat{X}_1^-$ is 
\begin{align}
	\widetilde{U}_{\hat{X}_1^-}=
	\begin{bmatrix}
		\frac{2}{2\sqrt{5}+\sqrt{2}} & 0 & * & \cdots & *\\  
		\frac{4}{2\sqrt{5}+\sqrt{2}} & 0 & * & \cdots & *\\  
		* & * & * & \cdots & * \\
		\vdots & \vdots & \vdots &  \ddots &  \vdots\\
		* & * & * & \cdots & *
	\end{bmatrix}\in BE_{2\sqrt{5}+\sqrt{2},3}(\hat{X}_1^-).
\end{align}

(3) In the similar way, the calculation of $\widetilde{U}_{A P_0 A^T}$ is done first, followed by obtaining the block encoding of $P_1^-$. The specific calculation process is shown as follows:
\begin{align}
	\widetilde{U}_{A P_0 A^T}=
	\begin{bmatrix}
		\frac{1}{2\sqrt{2}} & 0 & * & \cdots & *\\  
		0 & \frac{1}{2\sqrt{2}} & * & \cdots & *\\  
		* & * & * & \cdots & * \\
		\vdots & \vdots & \vdots &  \ddots &  \vdots\\
		* & * & * & \cdots & *
	\end{bmatrix}\in BE_{4\sqrt{2},3}(A P_0 A^T).
\end{align}
Let $T=4\sqrt{2}\widetilde{U}_{A P_0 A^T}+\sqrt{2}\widetilde{U}_Q$, the block encoding of $P_1^-$ is
\begin{align}
	\widetilde{U}_{P_1^-}=
	\begin{bmatrix}
		\frac{3}{5\sqrt{2}} & 0 & * & *\\  
		0 & \frac{3}{5\sqrt{2}} & * & *\\ 
		* & * & * &  * \\
		* & * & * &  *
	\end{bmatrix}\in BE_{5\sqrt{2},4}(P_1^-).	
\end{align}

(4) The calculation of the Kalman gain poses one of the primary challenges within the entire algorithm.

Firstly, we conduct matrix multiplication operations to obtain the block encodings of $P_1^-H^T$ and $H P_1^- H^T$, respectively.
\begin{equation}
	\begin{aligned}
		\widetilde{U}_{P_1^-H^T}=
		\begin{bmatrix}
			\frac{6}{5\sqrt{10}} & 0 & * &\cdots & *\\  
			0 & \frac{3}{5\sqrt{10}}  & * &\cdots & *\\ 
			* & * & * & \cdots & * \\
			\vdots & \vdots & \vdots &  \ddots &  \vdots\\
			* & * & * & \cdots & *
		\end{bmatrix}\in BE_{5\sqrt{10},5}(P_1^-H^T). \\		
		\widetilde{U}_{H P_1^- H^T}=
		\begin{bmatrix}
			\frac{12}{25\sqrt{2}} & 0 & * &\cdots & *\\  
			0 & \frac{3}{25\sqrt{2}} & * &\cdots & *\\ 
			* & * & * & \cdots & * \\
			\vdots & \vdots & \vdots &  \ddots &  \vdots\\
			* & * & * & \cdots & *
		\end{bmatrix}\in BE_{25\sqrt{2},6}(H P_1^- H^T).
	\end{aligned}
\end{equation}

Secondly, by setting $\widetilde{U}_{A_{temp}}=25\sqrt{2}\widetilde{U}_{H P_1^- H^T}+\sqrt{2}\widetilde{U}_R$, we can derive the block encoding for $A_{temp}$.
\begin{align}
	\widetilde{U}_{A_{temp}}=
	\begin{bmatrix}
		\frac{1}{2\sqrt{2}} & 0 & * &\cdots & *\\  
		0 & \frac{2}{13\sqrt{2}} & * &\cdots & *\\ 
		* & * & * & \cdots & * \\
		\vdots & \vdots & \vdots &  \ddots &  \vdots\\
		* & * & * & \cdots & *
	\end{bmatrix}\in BE_{26\sqrt{2},7}(A_{temp}).		
\end{align}
The specific circuit for $\widetilde{U}_{A_{temp}}$ is illustrated in Figure \ref{fig:Circuit_A_temp}. At last, $A_{temp}$ can be obtained by measuring the final state. The sampling process is 10 iterations, and 8192 shots are executed in each iteration. Of course, each sampling result has a small difference, and we take the average for the input of the subsequent process. The selected result of $A_{temp}$ is:

\begin{align}
	A_{temp}=
	\begin{bmatrix}
		13.0000 & 0 \\  
		0 & 4.0010 
	\end{bmatrix}
\end{align}

The value is obtained as 
\begin{align}
	A_{{temp},c}=
	\begin{bmatrix}
		13 & 0 \\  
		0 & 4 
	\end{bmatrix}
\end{align} by the classical algorithm. The solution closely aligns with the classical algorithm, showing only a 0.1\% discrepancy upon comparison. However, there are occasional fluctuations in the sampling that may result in an error of 0.01.

\begin{figure}[!h]
	\centering
	\includegraphics[scale=0.2]{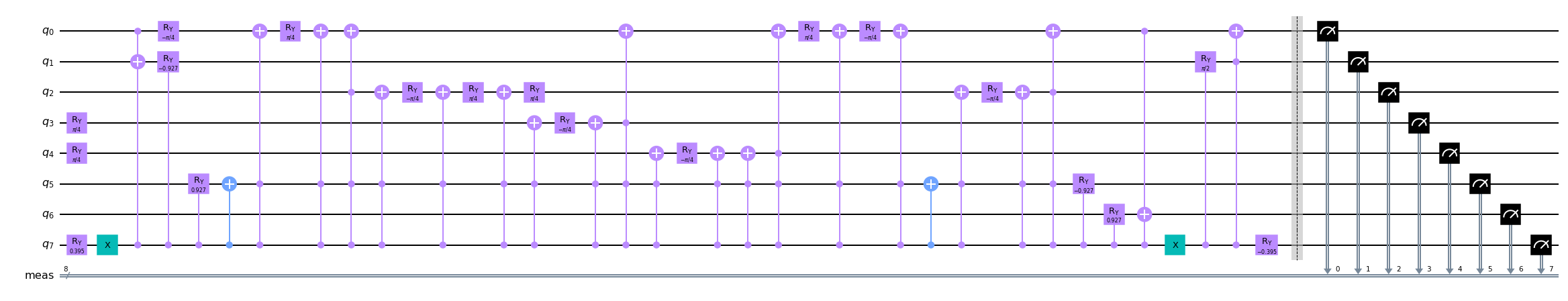}
	\renewcommand{\figurename}{Figure} 
	\caption{The quantum circuit of calculating $\widetilde{U}_{A_{temp}}$}
	\label{fig:Circuit_A_temp}
\end{figure} 
The third step involves calculating the matrix inversion of $A_{temp}$ using the QSVT method. It is worth noting that the condition number of $A_{temp}$ is $\kappa\approx3.25$. To accomplish this, we utilize pyqsp\cite{Chuang2022pyqsp}, a Python package specifically designed for generating QSP (quantum signal processing) phase angles. In our case, we set $\kappa=3.5$ and $\epsilon=0.01$. Consequently, we obtain a comprehensive set of QSP phase angles using Chebyshev polynomials that effectively approximate $\frac{1}{x}$ as follows:
\begin{equation}
		\begin{aligned}
		 \Phi=[-&0.06009860, \hspace{0.2em}-0.22167921,   \hspace{0.2em}-0.55045474,   \hspace{0.2em}-0.77609938, \hspace{0.2em}-0.52660613,\hspace{1em}0.21660955,\\
		  &0.38314397,  \hspace{1em} 0.56345849, \hspace{1em} 0.61656709, \hspace{0.2em}-0.40620100,\hspace{0.2em}-0.69935839, \hspace{1em}0.87923605, \\
		  &0.07902779,  \hspace{1em} 0.05022705, \hspace{1em} 1.11206017,\hspace{0.2em}-2.71868915, \hspace{1em}2.51891920,  \hspace{0.2em}-0.62276767,\\
		  -&0.52968031,  \hspace{0.2em} -0.11403588,  \hspace{0.2em}-0.72213012,   \hspace{0.2em}-0.03486492, \hspace{0.2em}-0.37990064, \hspace{1em}0.80318098, \\ 
				& 0.49191184,  \hspace{1em} 1.15950162, \hspace{0.2em}-0.60287972,  \hspace{0.2em}-0.60286045, \hspace{1em}1.15950755,\hspace{0.2em}-2.64969722,  \\
				& 0.80316442,  \hspace{0.2em}-0.37990169,   \hspace{0.2em}-0.03487192,   \hspace{1em}2.41946175,\hspace{0.2em}-0.11403621, \hspace{0.2em}-0.52967421,  \\
				-&0.62277180,  \hspace{0.2em}-0.62267939,   \hspace{1em}0.42290199,  \hspace{1em}1.11205726,\hspace{1em}0.05023048,\hspace{1em}0.07903506,  \\
				-&2.26235974,  \hspace{1em} 2.44222583, \hspace{0.2em}-0.40619834,   \hspace{0.2em}-2.52502044, \hspace{1em}0.56345461,\hspace{1em}0.38314771,  \\
				& 0.21661194,  \hspace{0.2em} -0.52660677,  \hspace{1em} 2.36549409, \hspace{0.2em}-0.55045294, \hspace{0.2em}-0.22168027, \hspace{1em}1.51069668].\nonumber
	\end{aligned} 
\end{equation}
denoted by $\Phi=[\phi_0,\phi_1,\cdots,\phi_{53}]$.

The corresponding polynomial approximation is shown in Figure \ref{fig:Inversion}. The range of x is $[-1,-\frac{1}{3.5}] \cup [\frac{1}{3.5},1]$, the scaling factor of $f(x)$ is 5.40127692.

\begin{figure}[!h]
	\centering
	\includegraphics[scale=0.3]{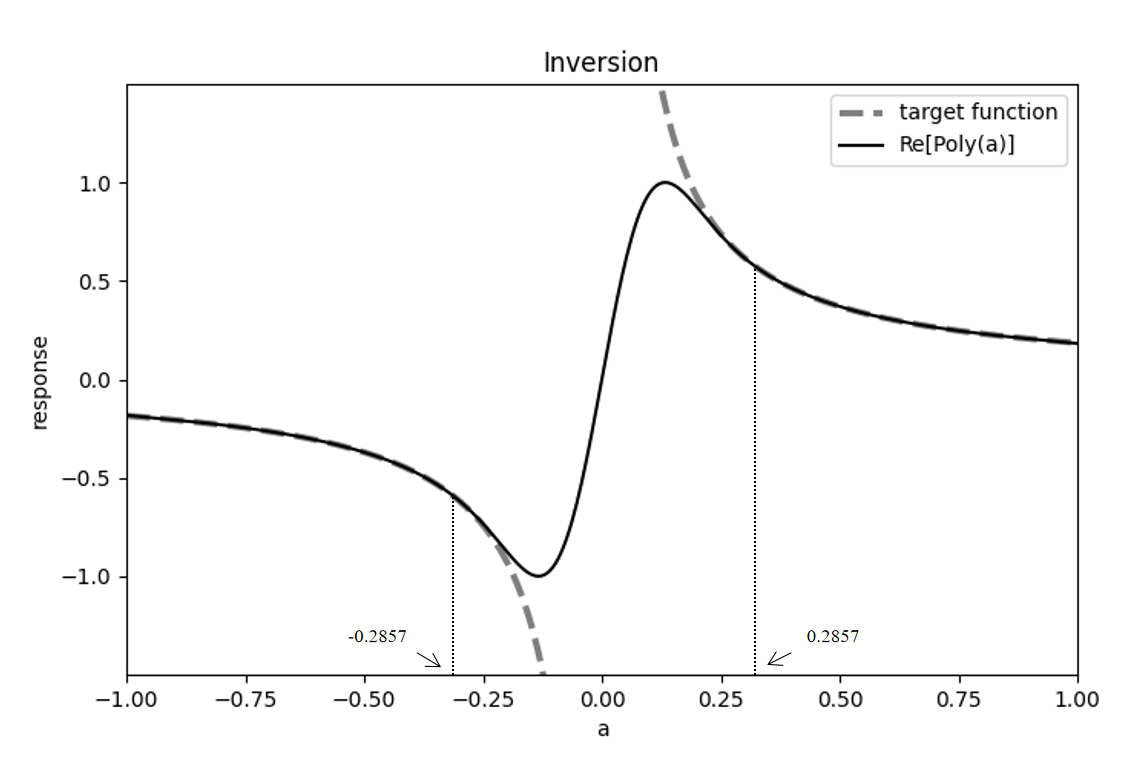}
	\renewcommand{\figurename}{Figure} 
	\caption{The function for the polynomial approximation to the inverse function with $\kappa = 3.5$ and $\epsilon = 0.01$}
	\label{fig:Inversion}
\end{figure} 

Apparently, the singular values of $A_{temp}$ are $\sigma_1=\frac{1}{2\sqrt{2}},\sigma_2=\frac{2}{13\sqrt{2}}$, but $\sigma_2$ is out of the range. We need to rescale $A_{temp}$ to obey this condition. The original matrix of $A_{temp}$ is 
$\begin{bmatrix}
	13 & 0 \\  
	0 & 4 
\end{bmatrix}$, we can reconstruct the block encoding of $A_{temp}$, denoted by $\widetilde{U}_{A_{temp}}$.

\begin{align}
	\widetilde{U}_{A_{temp}}=
	\begin{bmatrix}
		\frac{13}{\sqrt{185}} & 0 & * & *\\  
		0 & \frac{4}{\sqrt{185}} & *  & *\\ 
		* & * & * & * \\
		* & * & * & *
	\end{bmatrix}\in BE_{\sqrt{185},2}(A_{temp}).
\end{align}	
Obviously, the singular values of the new matrix $A_{temp}$ are evidently within the specified range.	

It's easy to get $U_{A_{temp}}^\dagger$ through $U_{A_{temp}}$. Then we can get the block encoding of $P(A_{temp})$ by constructing the circuit shown in Figure \ref{fig:Circuit_A_inv}. $d=53$, we construct the parameterized circuit consisting of $U_{A_{temp}},U_{A_{temp}}^\dagger$ are applied alternately by 53 iterates and 54 rotations that are interleaved in alternating sequence.
\begin{align}
	\widetilde{U}_{P(A_{temp})}=
	\begin{bmatrix}
		0.19306 & 0 & * &\cdots & *\\  
		0 & 0.62846 & * &\cdots & *\\ 
		* & * & * & \cdots & * \\
		\vdots & \vdots & \vdots &  \ddots &  \vdots\\
		* & * & * & \cdots & *
	\end{bmatrix}\in BE_{\frac{1}{\sqrt{185}}\times 5.40127692,2}(A_{temp}^{-1},0.001).		
\end{align}

\begin{figure}[!h]
	\centering
	\includegraphics[scale=0.35]{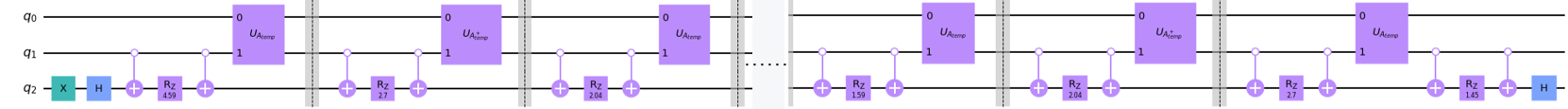}
	\renewcommand{\figurename}{Figure} 
	\caption{The quantum circuit of calculating $\widetilde{U}_{A_{temp}^{-1}}$}
	\label{fig:Circuit_A_inv}
\end{figure} 

Lastly, we calculate the block encoding of $K_1=P_1^-H^T A_{temp}^{-1}$.
\begin{align}
	\widetilde{U}_{K_1}=
	\begin{bmatrix}
		0.07326 & 0 & * &\cdots & *\\  
		0 & 0.11924 & * &\cdots & *\\ 
		* & * & * & \cdots & * \\
		\vdots & \vdots & \vdots &  \ddots &  \vdots\\
		* & * & * & \cdots & *
	\end{bmatrix}\in BE_{\frac{5\sqrt{74}}{37}\times 5.40127692,7}(K_1,0.001).		
\end{align}

(5) Calculate $\hat{X}_1=\hat{X}_1^-+K_1(Z_1-H\hat{X}_1^-)$. Firstly, we can get the block encoding of $H\hat{X}_1^-$,
\begin{align}
	\widetilde{U}_{H\hat{X}_1^-}=
	\begin{bmatrix}
		\frac{4}{10+\sqrt{10}} & 0 & * &\cdots & *\\  
		\frac{4}{10+\sqrt{10}} & 0 & * &\cdots & *\\ 
		* & * & * & \cdots & * \\
		\vdots & \vdots & \vdots &  \ddots &  \vdots\\
		* & * & * & \cdots & *
	\end{bmatrix}\in BE_{10+\sqrt{10},4}(H\hat{X}_1^-).	
\end{align}
Secondly, we perform the block encoding of $Z_1-H\hat{X}_1^-$,

\begin{align}
	\widetilde{U}_{Z_1-H\hat{X}_1^-}=
	\begin{bmatrix}
		-0.2056 & 0 & * &\cdots & *\\  
		-0.2071 & 0 & * &\cdots & *\\ 
		* & * & * & \cdots & * \\
		\vdots & \vdots & \vdots &  \ddots &  \vdots\\
		* & * & * & \cdots & *
	\end{bmatrix}\in BE_{10+\sqrt{10}+\sqrt{2},5}(Z_1-H\hat{X}_1^-).
\end{align}
And then we multiply $\widetilde{U}_{k_1}$ times $\widetilde{U}_{Z_1-H\hat{X}_1^-}$ to get the block encoding of $K_1(Z_1-H\hat{X}_1^-)$.
\begin{align}
	\widetilde{U}_{K_1(Z_1-H\hat{X}_1^-)}=
	\begin{bmatrix}
		-0.01507 & 0 & * &\cdots & *\\  
		-0.02454 & 0 & * &\cdots & *\\ 
		* & * & * & \cdots & * \\
		\vdots & \vdots & \vdots &  \ddots &  \vdots\\
		* & * & * & \cdots & *
	\end{bmatrix}\in BE_{91.5237,12}(K_1(Z_1-H\hat{X}_1^-),0.001).
\end{align}
At last, perform matrix addition to get the block encoding of $\hat{X}_1$. Let $T=(2\sqrt{5}+\sqrt{2})\widetilde{U}_{\hat{X}_1^-}+91.5237 \\
\widetilde{U}_{K_1(Z_1-H\hat{X}_1^-)}$, the result is
\begin{align}
	\widetilde{U}_{\hat{X}_1}=
	\begin{bmatrix}
		0.00636 & 0 & * &\cdots & *\\  
		0.01801 & 0 & * &\cdots & *\\ 
		* & * & * & \cdots & * \\
		\vdots & \vdots & \vdots &  \ddots &  \vdots\\
		* & * & * & \cdots & *
	\end{bmatrix}\in BE_{97.41,13}(\hat{X}_1,0.001).	
\end{align}	

\begin{figure}[!h]
	\centering
	\includegraphics[scale=0.25]{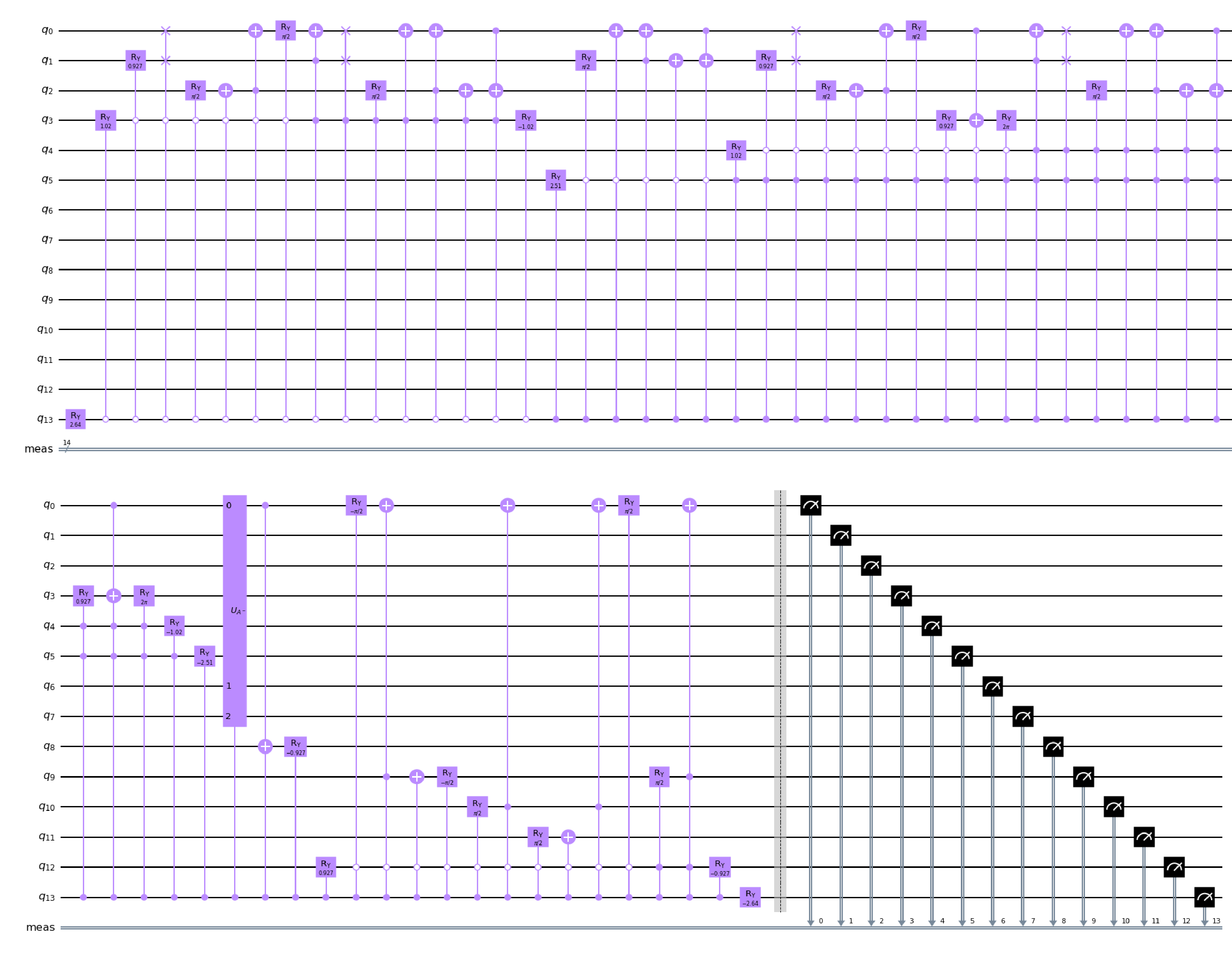}
	\renewcommand{\figurename}{Figure} 
	\caption{The quantum circuit of calculating $\widetilde{U}_{\hat{X}_1}$}
	\label{fig:Circuit_X_k}
\end{figure}

(6) Calculate $P_1=P_1^--K_1HP_1^-$. Referring to the previous process, firstly, we can construct the circuit to get the block encoding of $-K_1HP_1^-$.
\begin{align}
	\widetilde{U}_{-K_1HP_1^-}=
	\begin{bmatrix}
		-0.0278 & 0 & * &\cdots & *\\  
		0 & -0.02262 & * &\cdots & *\\ 
		* & * & * & \cdots & * \\
		\vdots & \vdots & \vdots &  \ddots &  \vdots\\
		* & * & * & \cdots & *
	\end{bmatrix}\in BE_{99.2774,12}(-K_1HP_1^-,0.001).	
\end{align}
And then we add $\widetilde{U}_{-K_1HP_1^-}$ to $\widetilde{U}_{P_1^-}$ to get the block encoding of $P_1$.
\begin{align}
	\widetilde{U}_{P_1}=
	\begin{bmatrix}
		2.2567\times 10^{-3} & 0 & * &\cdots & *\\  
		0 & 7.0889\times 10^{-3} & * &\cdots & *\\ 
		* & * & * & \cdots & * \\
		\vdots & \vdots & \vdots &  \ddots &  \vdots\\
		* & * & * & \cdots & *
	\end{bmatrix}\in BE_{106.3485,13}(P_1,0.001).	
\end{align}

\begin{figure}[!h]
	\centering
	\includegraphics[scale=0.25]{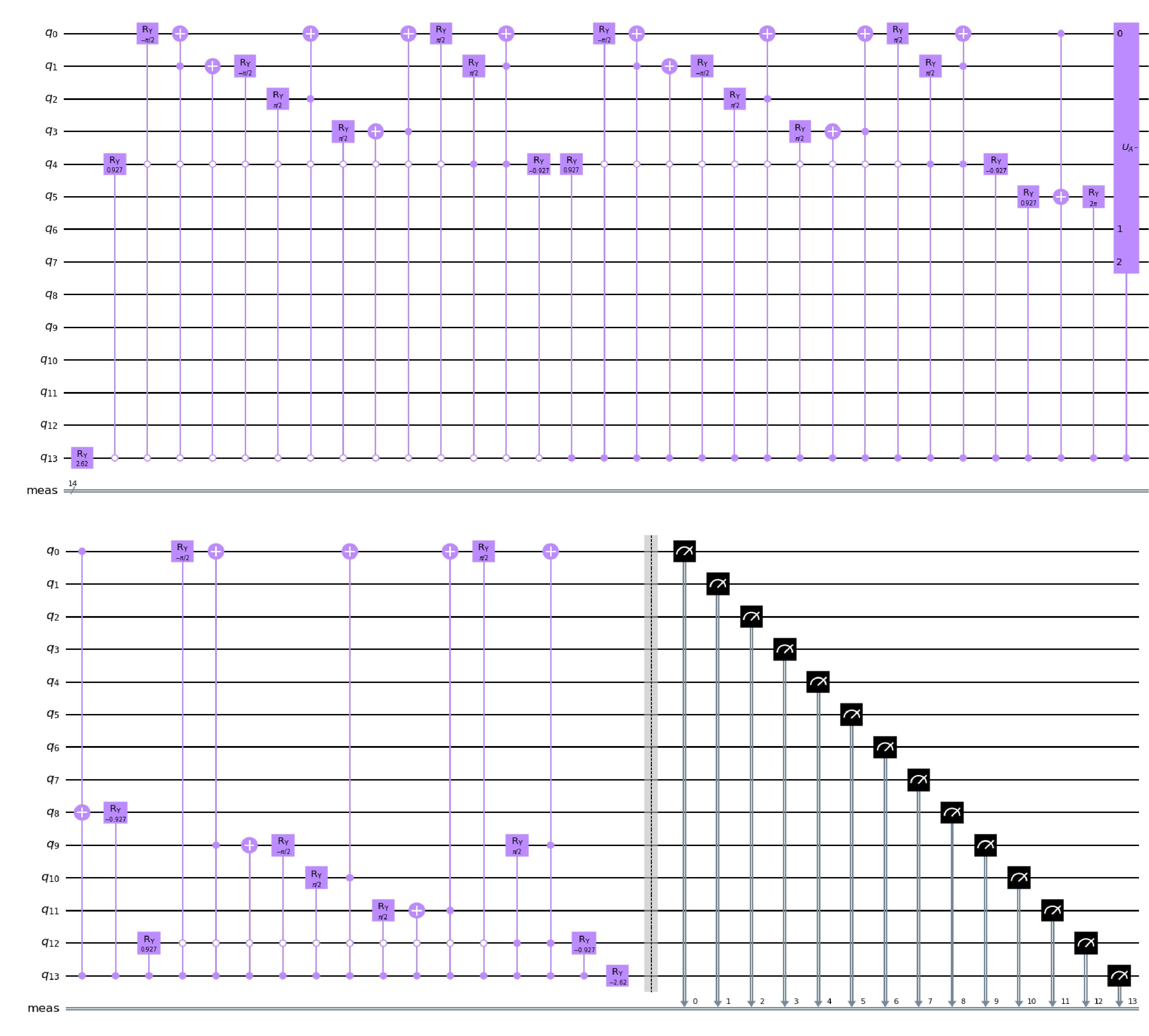}
	\renewcommand{\figurename}{Figure} 
	\caption{The quantum circuit of calculating $\widetilde{U}_{P_1}$}
	\label{fig:Circuit_P_k}
\end{figure}

\section{Experimental results and discussion}
In order to obtain $\hat{X}_1$ and $P_1$, it is essential to perform measurements on $\widetilde{U}_{\hat{X}_1}$ and $\widetilde{U}_{P_1}$. To compare the estimated amplitudes with the exact ones, it is crucial to represent the final state in its full-amplitude form. 

\subsection{The form of the final state}

For $\widetilde{U}_{\hat{X}_1}$, we set the initial state to be $|0^{\otimes 14}\rangle$, after $\widetilde{U}_{\hat{X}_1}$ acting on $|0^{\otimes 14}\rangle$ the final state $|\psi_{\hat{X}_1} \rangle$ is:
\begin{equation}
	\begin{aligned}
		|\psi_{\hat{X}_1} \rangle&=6.36484546\times10^{-3}|00\cdots00\rangle+1.80050738\times10^{-2}|00\cdots01\rangle \\
		&-8.95345618\times10^{-16}|00\cdots10\rangle+\cdots+3.17131363\times10^{-19}|11\cdots01\rangle\\
		&+4.02491466\times10^{-24}|11\cdots10\rangle+1.10222235\times10^{-23}|11\cdots11\rangle.
	\end{aligned}
\end{equation}

Similarly, for $\widetilde{U}_{P_1}$, the initial states are set to $|0^{\otimes 14}\rangle,|0^{\otimes 13}\rangle|1\rangle$ respectively, after the operation of $\widetilde{U}_{P_1}$ the two final states are respectively:
\begin{equation}
	\begin{aligned}
		|\psi_{P_{11}} \rangle&=2.25672976\times 10^{-3}|00\cdots00\rangle+2.10564344\times 10^{-15}|00\cdots01\rangle \\
		&+1.50448651\times10^{-3}|00\cdots10\rangle+\cdots-6.73134796\times10^{-19}|11\cdots01\rangle \\
		&-1.38326258\times10^{-21}|11\cdots10\rangle-8.18294411\times10^{-19}|11\cdots11\rangle.
	\end{aligned}
\end{equation}
\begin{equation}
	\begin{aligned}
		|\psi_{P_{12}} \rangle&=3.14223305\times 10^{-16}|00\cdots00\rangle+7.08894955\times 10^{-3}|00\cdots01\rangle \\
		&+7.52233295\times10^{-4}|00\cdots10\rangle+\cdots+1.61401008\times10^{-20}|11\cdots01\rangle \\
		&+1.85712631\times10^{-18}|11\cdots10\rangle+2.03004883\times10^{-20}|11\cdots11\rangle.
	\end{aligned}
\end{equation}

\subsection{Measurement}

$|\psi_X \rangle$ consists of $2^{14}=16384$ quantum superposition states, the experiment involves conducting $16384$ shots in one sampling iteration. We need 100 iterations for an acceptable result. The resulting probabilities of partial samples in $|\psi_{\hat{X}_1} \rangle$ are shown in Figure \ref{fig:Probability_X_k}. The two values of $\hat{X_1}$ can be obtained by squaring $P(|00\cdots00\rangle)$ and $P(|00\cdots01\rangle)$ and multiplying by the normalization factor 97.41:
\begin{align}
	\hat{X}_1=
	\begin{bmatrix}
		0.6201 \\  
		1.7635
	\end{bmatrix}.
\end{align}

The probabilities of partial samples in $|\psi_{P_{11}}\rangle$ are illustrated in Figure \ref{fig:Probability_P1}. The first row of $P_1$ can be obtained by squaring the probabilities for $|00\cdots00\rangle$ and $|00\cdots01\rangle$ and multiplying the normalization factor 106.3485, respectively. Additionally, Figure \ref{fig:Probability_P2} displays the probabilities of partial samples in $|\psi_{P_{12}}\rangle$. Similarly, the second row of $P_1$ can be derived by the same way:
\begin{align}
	P_1=
	\begin{bmatrix}
		0.2213 & 0 \\ 
		0 & 0.7619
	\end{bmatrix}.
\end{align}

\begin{figure}[!t]
	\centering
	\includegraphics[scale=0.2]{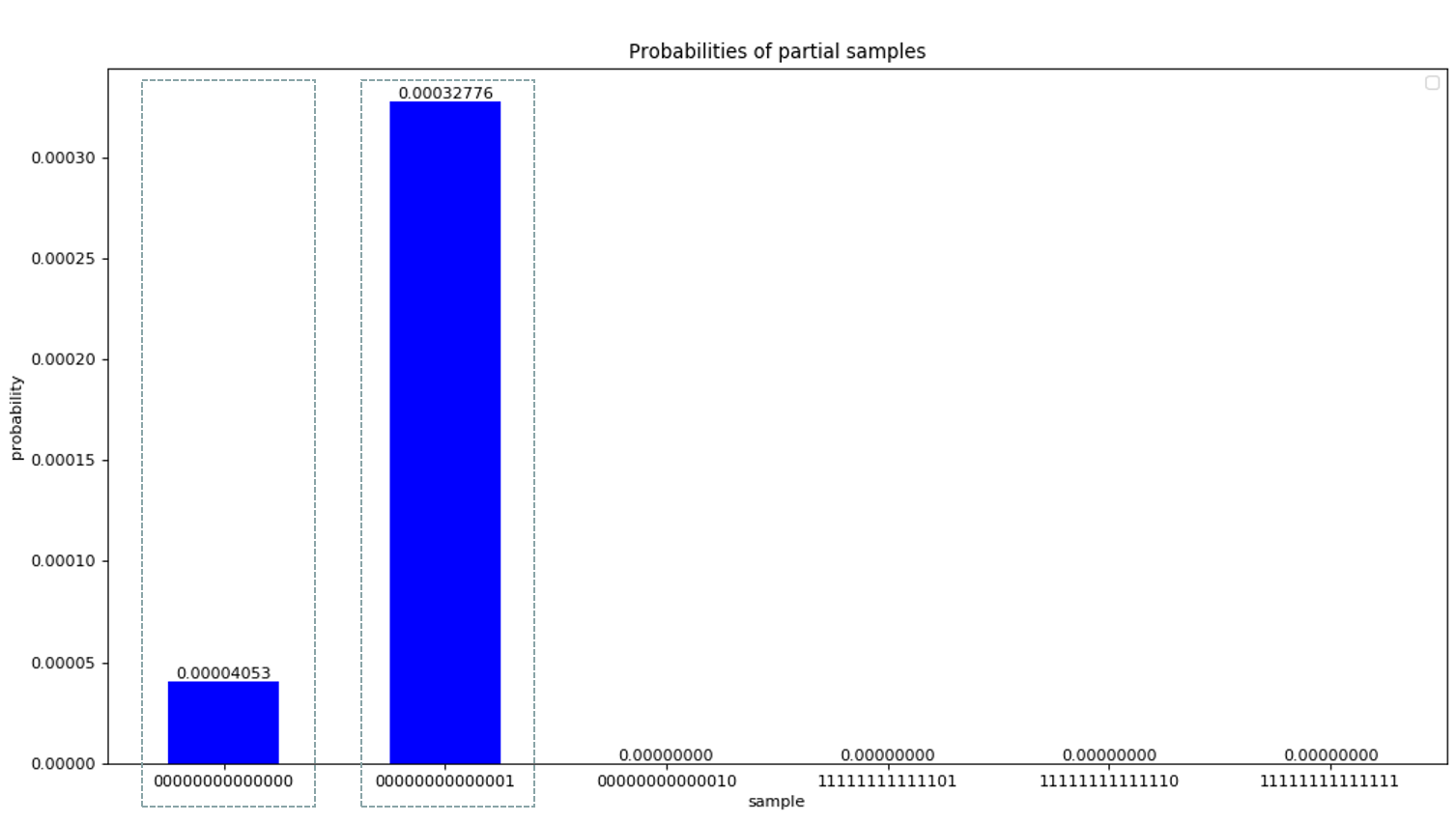}
	\renewcommand{\figurename}{Figure} 
	\caption{The probabilities of partial samples in $|\psi_{\hat{X}_1} \rangle$}
	\label{fig:Probability_X_k}
\end{figure} 

\begin{figure}[!t]
	\centering
	\begin{minipage}{0.45\linewidth}
		\includegraphics[width=\linewidth]{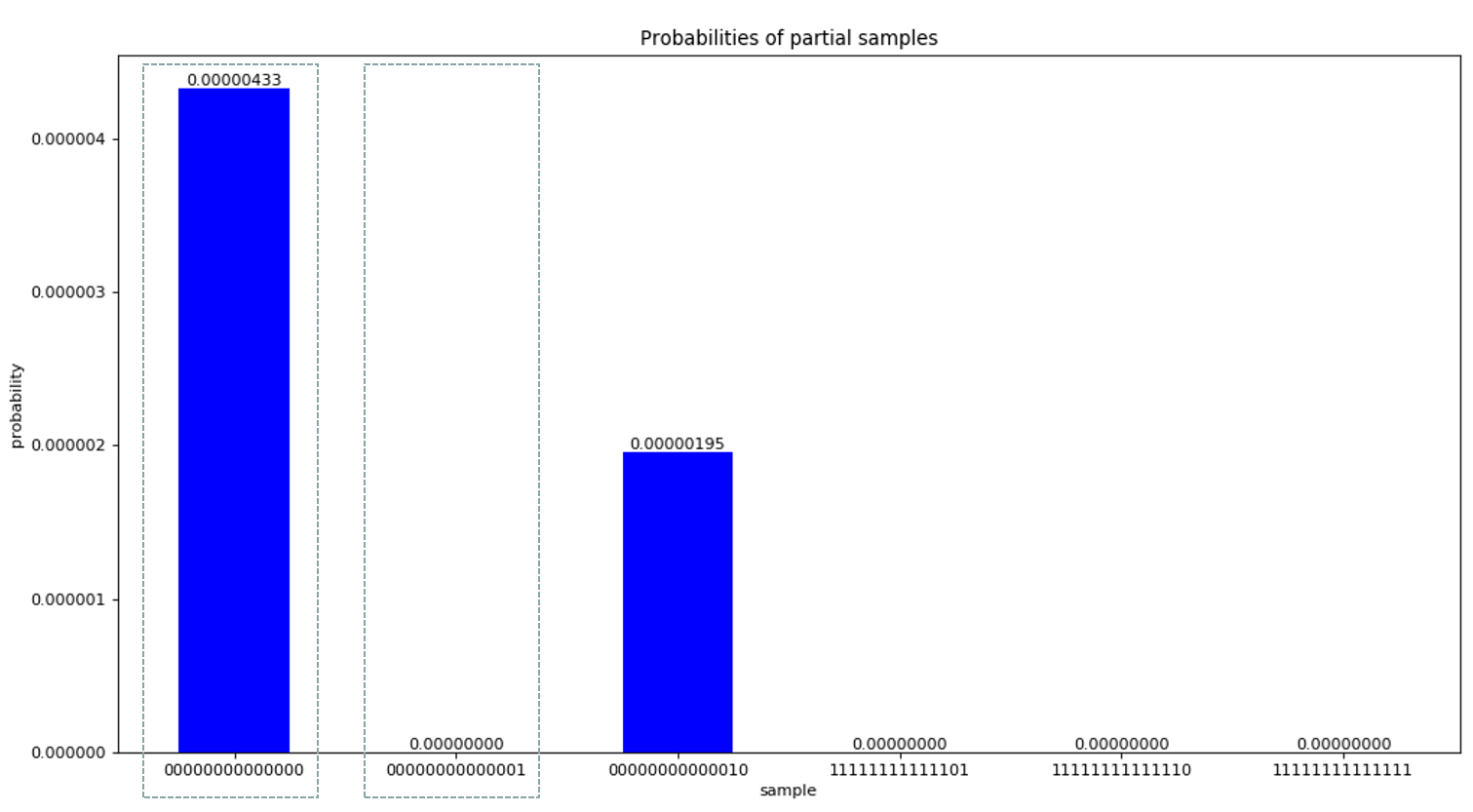}
		\renewcommand{\figurename}{Figure} 
		\caption{The probabilities of partial samples in $|\psi_{P_{11}} \rangle$}
		\label{fig:Probability_P1}
	\end{minipage}
	\begin{minipage}{0.45\linewidth}
		\includegraphics[width=\linewidth]{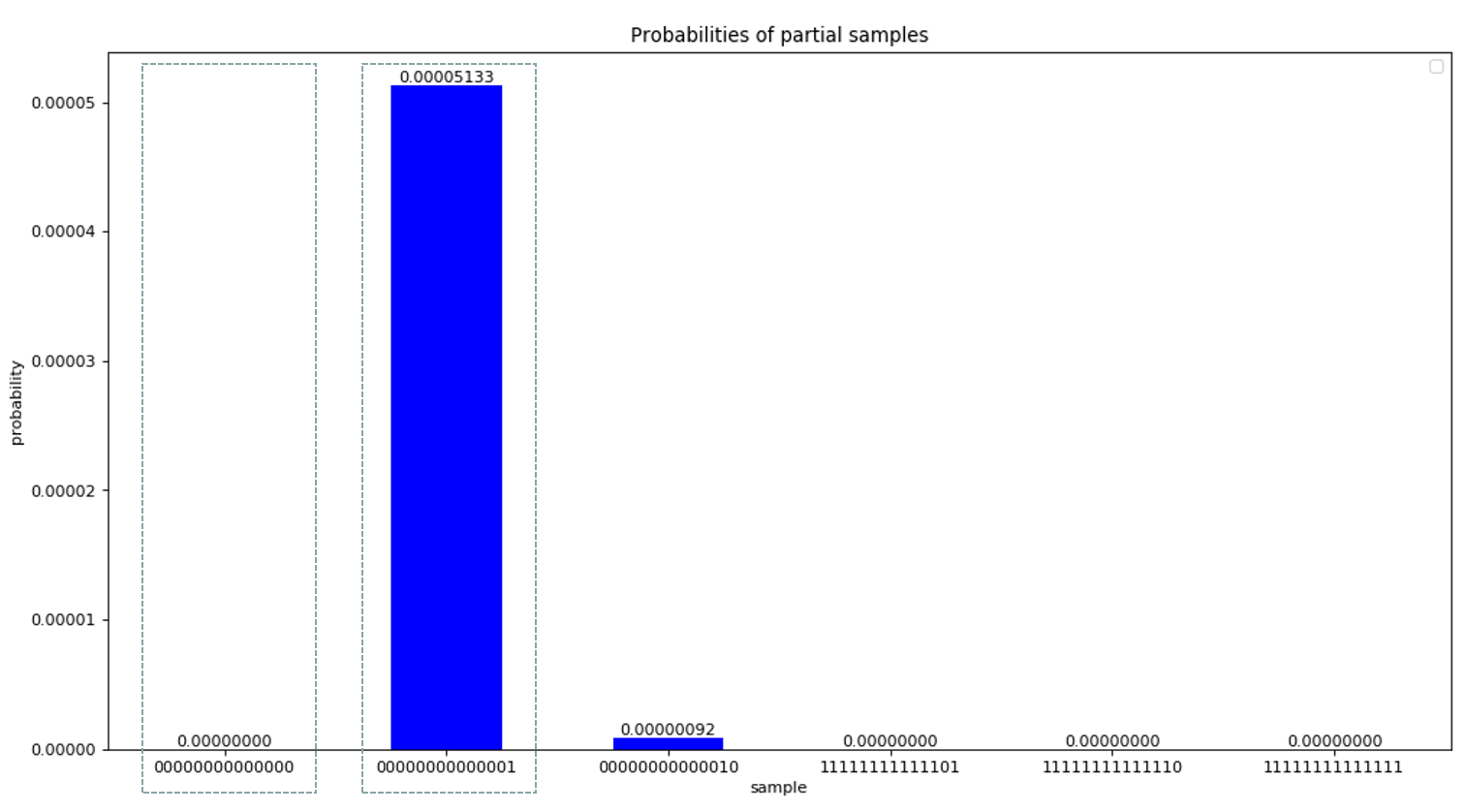}
		\renewcommand{\figurename}{Figure} 
		\caption{The probabilities of partial samples in $|\psi_{P_{12}} \rangle$}
		\label{fig:Probability_P2}
	\end{minipage}
	
\end{figure}

The results of $\hat{X}_1$ and $P_1$ obtained through classical approach are:
\begin{align}
	\hat{X}_{1,c}=
	\begin{bmatrix}
		0.6154 \\  
		1.75
	\end{bmatrix},
	P_{1,c}=
	\begin{bmatrix}
		0.2308 & 0 \\ 
		0 & 0.75
	\end{bmatrix}.
\end{align}
Due to the random nature of sampling, errors are inevitable in the results obtained from 100 iterations. In our study, we implemented 10 sets of experiments with each set consisting of 100 and 10,000 iterations respectively (shots=16384). A comparative analysis between classical and quantum approaches was performed, and the results are illustrated in Figure \ref{fig:Comparation_X1} and \ref{fig:Comparation_P1}. Here $\hat{x}_0$ and $\hat{x}_1$ represent the two elements of vector $\hat{X}_1$, while $P_{11}$ and $P_{22}$ denote two non-zero entries in matrix $P_1$. The outcomes obtained from 10000 iterations closely resemble classical results with a negligible error margin of only 0.01. Conversely, values derived from 100 iterations exhibit a larger error margin of 0.1.

\begin{figure}[!t]
	\centering
	\begin{minipage}{0.45\linewidth}
		\includegraphics[width=\linewidth]{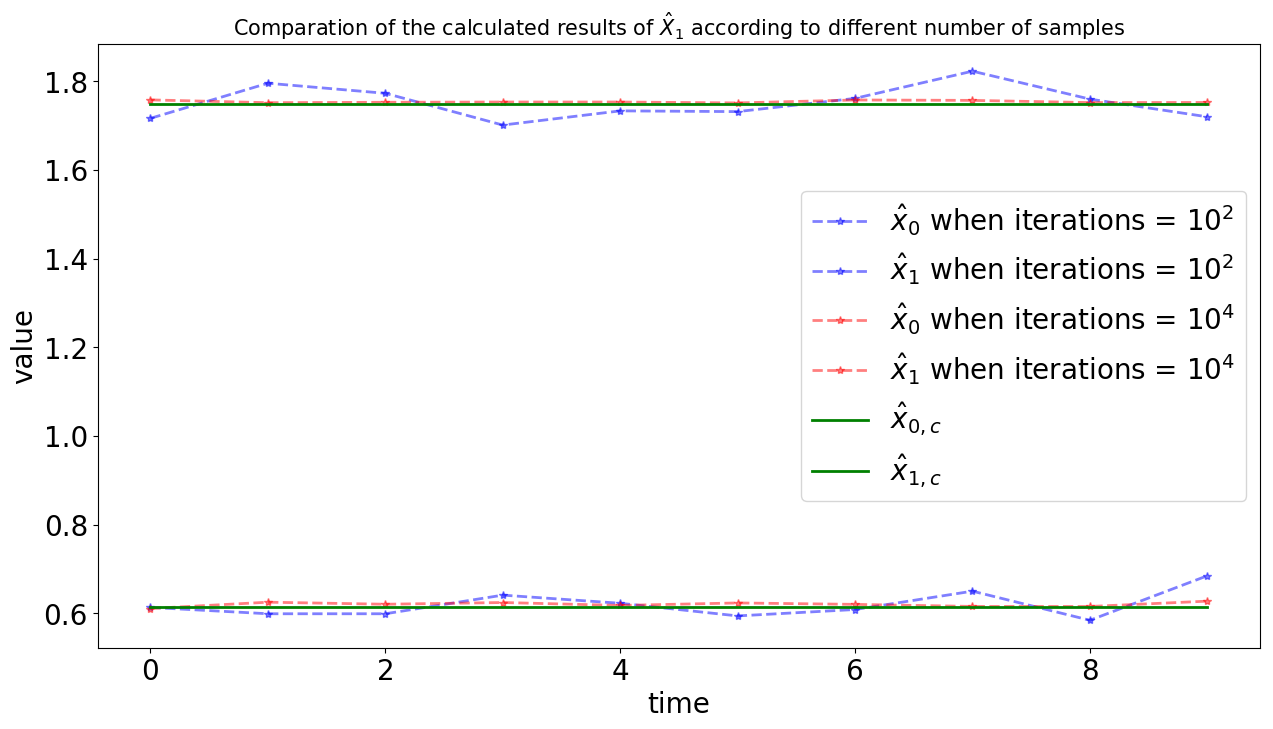}
		\renewcommand{\figurename}{Figure} 
		\caption{Comparation of $\hat{X}_1$ between quantum and classical approaches}
		\label{fig:Comparation_X1}
	\end{minipage}
	\begin{minipage}{0.45\linewidth}
		\includegraphics[width=\linewidth]{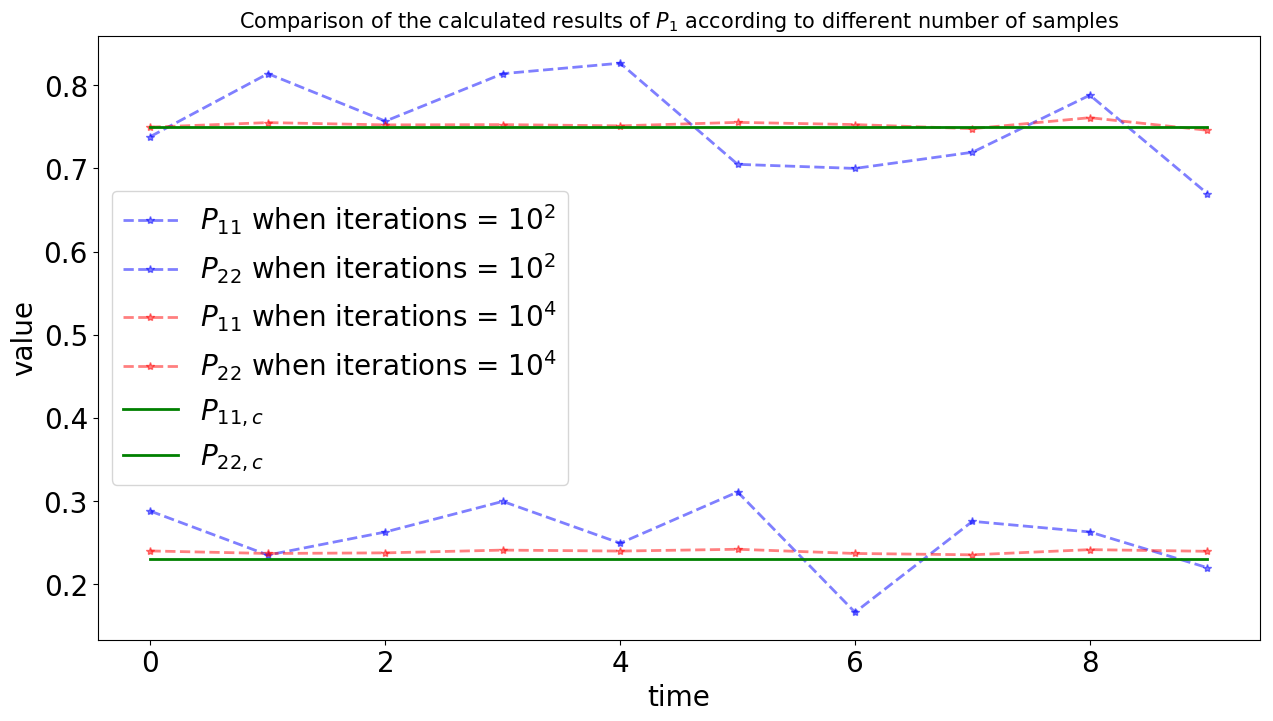}
		\renewcommand{\figurename}{Figure} 
		\caption{Comparation of $P_1$ between quantum and classical approaches}
		\label{fig:Comparation_P1}
	\end{minipage}	
\end{figure}

\subsection{Error analysis}

The previous subsection presents an illustrative example showcasing a quantum verison of the Kalman filter algorithm using block encoding. By comparing $\hat{X_1}, P_1$ obtained from quantum and classical algorithms, we observe that the quantum algorithm achieves an error. This error comes from three aspects when the block encoding is completely accurate: 

(1)The algorithm itself. For $\hat{X_1}$, we can directly calculate the two elements from the amplitude of $|0^{\otimes14}\rangle,|0^{\otimes13}\rangle|1\rangle$, which are 0.6199 and 1.7519 respectively. The error is 0.001 compared with the result obtained by the classical approach. A similar conclusion holds for $P_1$. The error from the algorithm itself mainly coming from the results of matrix inversion leading to an error in Kalman gain $K_k$. This is unavoidable because of the polynomial approximation, and we can only increase the degrees of freedom $d$ to find a Chebyshev polynomial that more closely approximates $f(x)=1/x$ and reduce the error as much as possible. In many cases, the matrix inversion also exist errors in the classical algorithm.  

(2) Number of sampling. The errors of the results obtained with 100 and 10000 iterations are different. We can mitigate the error of the second kind to increase the number of iterations or conduct multiple sets of experiments. However, the number of samples does not increase without bound either. The accuracy is directly proportional to the square root of the sample size, as it is well-established in scientific literature. The quantum state $|\psi_X\rangle,|\psi_{P_1}\rangle,|\psi_{P_2}\rangle$ consists of 14 qubits, with a total of $2^{14}$ superposition states. Obtaining accurate probability values through sampling at such a precision level of $10^{-4}$ is challenging and would require approximately $10^8$ samples. But for solving general $2\times2$ matrices, this sampling scale is not necessary.

(3) Normalization factor. $\hat{X_0}$ is a $2$-dimensional column vector, but the dimension of $\widetilde{U}_{\hat{X}_1}$ becomes $2^{14}\times 2^{14}$, expanding by a factor of 13. This leads to a significantly large normalization factor for the block-encoded matrix $\widetilde{U}_{\hat{X}_1}$, reaching around 100 in this instance. We know that the real matrix equals the block-encoded matrix multiplied by the normalization factor. However, even with an error of only 0.0001 in the block-encoded matrix, it translates to an error of 0.01 in the real matrix, which is unavoidable.

\subsection{Discussion}
Furthermore, it is important to consider certain caveats associated with this approach.

(1) The accuracy improvement achieved with $10^8$ samples is only 0.01, indicating compromised precision for faster computation in this methodology. When comparing the time complexity of quantum and classical algorithms, it's often overlooked that both the encoding and measurement processes are time-consuming. Some quantum algorithms require significant time resources for encoding and measurement procedures. In this algorithm, we have considered the temporal requirements for block encoding in QRAM while maintaining a lower time complexity compared to classical algorithms. However, if we also consider the process of measurements, it may not necessarily surpass classical algorithms as achieving more accurate results with this algorithm requires extensive measurements.

(2) The accumulation of errors. In the previous subsection, we compared the error in calculating $\hat{X}_{k-1}$ and $P_{k-1}$ between the classical and quantum algorithms, which yielded small differences. Furthermore, since $\hat{X}_{k-1}$ and $P_{k-1}$ serve as inputs for the subsequent loop, this implies that the outcomes of both algorithms will diverge. However, it is difficult to assess the impact of the exact solutions from the classical algorithm and approximate solutions from the quantum algorithm on the convergence of the Kalman filter. In practice, convergence is influenced by factors such as system model accuracy, data quality, and parameter selection making it a complex system engineering. Comparing results between these two algorithms holds limited value directly. 

(3) When solving the Kalman filtering algorithm, a crucial task is constructing block encoding, which can be quite challenging to implement. There is a significant amount of work involved in creating block encoding for matrices, and typically specific encoding methods exist for certain types of matrices\cite{Clader2022,Wan2019}. While there are some generalized encoding methods available for sparse matrices, block encoding remains a difficult task for more complex matrices. 

\section{Conclusions}
In this paper, we proposed a quantum scheme to solve state estimation in the Kalman filter. By block encoding, we encoded matrices stored in the quantum-accessible data structure into unitary matrices. Then the quantum circuit of solving the state estimation of the Kalman filter can be constructed through addition and multiplication rules of block encoding. Finally, the optimal state estimate is obtained by measuring the quantum state. And our scheme can reduce the time complexity from $O(n^3)$ to $O(\kappa\log(\kappa/\epsilon')poly\log(n/\epsilon))$ without considering the measurement, where $n$ is matrix dimension, $\kappa$ is condition number, $\epsilon$ is the desired precision in block encoding, $\epsilon'$ is the desired precision in matrix inversion. That means that an exponential speed-ups can be realized in the quantum algorithm in comparison to the existing algorithm on classical computers. And of course, the time complexity of block encoding generally takes $O(polylog(n/\epsilon))$. If a matrix requires $s$ qubits, the total number of needed qubits is linearly dependent on $s$ in the whole algorithm. So the space complexity of the whole algorithm is $O(s)$. 

The method presented in this paper can be applied to solve state estimation of the Kalman filter. Meanwhile, we show the potential of quantum computing for large-scale control systems. However, there are several problems that need to be solved in the future. Firstly, we have verified the feasibility of this algorithm theoretically, and the experiment was performed in a quantum simulator. Due to the large number of qubits and the deep circuit, the verification cannot be performed by existing quantum computers. Secondly, the construction of block encoding are still challenging, and it takes a lot of resources to implement. If there is a fast block encoding mechanism,  the power of the method based on block encoding model will be greater. Furthermore, in quantum algorithms, block encoding is a good way to embed a matrix into a unitary matrix that a quantum computer can perform. We can use this method to solve other control problems, and thus expand the application field of quantum computing.

\section*{Acknowledgments}
This work was supported by the National Natural Science Foundation of China (No.61673389).

\bibliographystyle{unsrt}  
\bibliography{references}

\end{document}